\documentclass[11pt]{article} 
\usepackage{a4wide}
\usepackage{amsmath,amsfonts}

\newtheorem{theorem}{Theorem}

\newtheorem{proposition}{Proposition}
\newtheorem{lemma}{Lemma} 
\newtheorem{claim}{Claim}

\numberwithin{lemma}{section}
\numberwithin{claim}{section}
\numberwithin{equation}{section}
 
\begin{document}

\title{Asymptotic stability of solitons of the gKdV equations \\
with general nonlinearity
 \footnote{
This research was supported in part by the Agence Nationale de la Recherche (ANR ONDENONLIN).
}}
\author{Yvan Martel$^{(1)}$ 
  and  Frank Merle$^{(2)}$
  }
\date{\quad \\ (1)
 Universit\'e de Versailles Saint-Quentin-en-Yvelines,
Math\'ematiques \\
 45, av. des Etats-Unis,
 78035 Versailles cedex, France\\  
 martel@math.uvsq.fr\\
\quad  \\
(2)
Universit\'e de Cergy-Pontoise, IHES and CNRS, Math\'ematiques   \\
2, av. Adolphe Chauvin,
95302 Cergy-Pontoise cedex, France \\
 Frank.Merle@math.u-cergy.fr\\
}
\maketitle
 
\begin{abstract}

We consider the generalized Korteweg-de Vries equation
\begin{equation}\label{kdvabs}
    \partial_t u + \partial_x (\partial_x^2 u + f(u))=0, \quad
    (t,x)\in [0,T)\times \mathbb{R},
\end{equation}
 with general $C^3$ nonlinearity $f$. 
Under an explicit condition on $f$ and $c>0$,
there exists a solution  in the energy space $H^1$ of \eqref{kdvabs} of the  type $u(t,x)=Q_c(x-x_0-ct)$,  called soliton.

In this paper, under general assumptions on $f$ and $Q_c$,
we prove that the family of soliton solutions around $Q_c$
is asymptotically stable in some local sense in $H^1$, i.e.
if $u(t)$ is close to $Q_{c}$ (for all $t\geq 0$), 
then $u(t)$ locally  converges in the energy space
to some $Q_{c_+}$  as $t\to +\infty$. Note in particular
that we do not assume the stability of $Q_{c}$.
This result is  based on a rigidity property of equation
\eqref{kdvabs} around $Q_{c}$ in the energy space whose proof relies on the introduction
of a dual problem.
These results extend the main results in \cite{MMjmpa}, \cite{MMarchives},
\cite{MMnonlinearity} and \cite{yvanSIAM}, devoted to the pure power case.
\end{abstract}

\section{Introduction}
We consider the generalized Korteweg-de Vries (gKdV) equations:
\begin{equation}\label{fkdv}
    \partial_t u + \partial_x (\partial_x^2 u + f(u))=0, \quad
    (t,x)\in [0,T)\times \mathbb{R},
\end{equation}
for $u(0)=u_0\in H^1(\mathbb{R})$,
with  a  general $C^3$ nonlinearity $f$.
We assume that  there exists an integer $p\geq 2$   such that
\begin{equation}\label{surfas1}
f(u)= a u^p+f_1(u)\quad \text{where  $a>0$ and $\lim_{u\to 0}\left|\frac{f_1(u)}{u^p}\right|=0$.}
\end{equation}
This is the only assumption on $f$ in this paper.  Denote $F(s)=\int_0^s f(s') ds'$. 
 
Note that the local  Cauchy problem is well-posed
in $H^1$, using the arguments of Kenig, Ponce and Vega \cite{KPV89}, \cite{KPV}, see  Remark 3 below. 
Moreover, the following conservation laws holds for $H^1$ solutions:
\begin{align*}
& \int u^2(t)=\int u_0^2,\quad 
E(u(t))=\frac 12 \int (\partial_x u(t))^2 -  \int F(u(t)) = \frac 12 \int (\partial_x u_{0})^2 -  \int F(u_0)
\end{align*}

Recall that if $Q_c$ is a solution of
\begin{equation}\label{ellipticas1} 
Q_c''+f(Q_c)=cQ_c, \quad x\in \mathbb{R}, \quad Q_c\in H^1(\mathbb{R}),
\end{equation}
then $R_{c,x_0}(t,x)=Q_c(x-x_0-ct)$
is solution of \eqref{fkdv}. We call soliton such nontrivial traveling wave solution
of \eqref{fkdv}. 

By well-known results on equation \eqref{ellipticas1} (see section 2), there exists
$c_*(f)>0$ (possibly, $c_*(f)=+\infty$) defined by 
$$c_*(f)=\sup\{c>0 \text{ such that }
\forall c'\in (0,c), \text{ $\exists$ $Q_{c'}$ positive solution of 
\eqref{ellipticas1}}\}.$$
See Section 2 for another characterization of $c_*(f)$ related to $f$.

Recall that if a solution $Q_c>0$ of \eqref{ellipticas1} 
exists then $Q_c$ is the unique (up to translation)
positive solution of \eqref{ellipticas1}
and  can be chosen  even on $\mathbb{R}$ and  decreasing on $\mathbb{R^+}$.

\smallskip

The main result of this paper is the following:

 \begin{theorem}[Asymptotic stability]\label{TH1as1}
 Assume that $f$ is $C^3$ and satisfies \eqref{surfas1}.
Let $0<c_0<c_*(f)$.
There exists $\alpha_0>0$ such that if
 $u(t)$ is a global  $(t\geq 0)$ $H^1$ solution of \eqref{fkdv}
satisfying
\begin{equation}\label{th1bas1}
	\forall t\geq 0,\quad
	\inf_{r\in \mathbb{R}}\|u(t,.+ r)-Q_{c_0}\|_{H^1} < \alpha_0,
	\end{equation}
then the following hold.
\begin{enumerate}
\item Asymptotic stability in the energy space.
There exist  $t\mapsto c(t)\in (0,c_*(f))$, $t\mapsto \rho(t)\in \mathbb{R}$ such that
\begin{equation}
\label{th1-1}
u(t)-Q_{c(t)}(.-\rho(t))\to 0\quad \text{in $H^1(x>\tfrac {c_0}{10} t)$ as $t\to +\infty$.}
\end{equation}
\item Convergence of the scaling parameter.
Assume further that there exists $\sigma_0>0$ such that
$c\mapsto \int Q_c^2$ is not constant in any interval
$I\subset [c_0-\sigma_0,c_0+\sigma_0]$. 
Then, by possibly taking a smaller $\alpha_0>0$, 
there exits $c_+\in (0,c_*(f))$ such that $c(t)\to c_+$ as $t\to +\infty$.
\end{enumerate}
\end{theorem}

The main ingredient of the proof of Theorem \ref{TH1as1} is a rigidity theorem around
$Q_{c_0}$.

 \begin{theorem}[Nonlinear Liouville Property around $Q_{c_0}$]\label{PROP2as1}
Assume that $f$ is $C^3$ and satisfies \eqref{surfas1}. Let $0<c_0<c_*(f)$.
There exists $\alpha_0>0$ such that if $u(t)$ is a global $(t\in \mathbb{R})$ $H^1$ solution of \eqref{fkdv}
satisfying,  for some function $t\mapsto \rho(t)$, $C,\sigma>0$,
\begin{align}
&	\forall t\in \mathbb{R},\quad
		\|u(t,.+ \rho(t))-Q_{c_0}\|_{H^1} \leq \alpha_0, \label{stableas1}\\
&	\forall t,x \in \mathbb{R} ,\quad
		|u(t,x+ \rho(t))|\leq  C e^{- \sigma |x|},
			\label{decayas1}
\end{align}
then there exists $c_1>0$, $x_1\in \mathbb{R}$, such that
$$ 	\forall t,x \in \mathbb{R} ,\quad
		u(t,x)=Q_{c_1}(x-x_1-c_1t).
$$
\end{theorem}
Theorem 1 above  is fundamental in proving the main results of \cite{MMcol2},
concerning the problem of collision of two solitary waves for general 
KdV equations. 
Indeed, as a corollary of the proofs  of Theorem \ref{TH1as1}, Theorem \ref{PROP2as1}
and \cite{MMT}, we obtain asymptotic stability of multi-solitons, see Section 5 for  a precise result and more details on the proofs.
See also \cite{MMas2} for more qualitative properties.

The arguments of \cite{MMcol1} and \cite{MMcol2} allow to describe the collision of two
solitary waves in a large but fixed interval of time. Large time asymptotics are necessary to preserve the soliton structure after the collision as $t\to +\infty$
 (Theorem \ref{TH1as1} and Theorem \ref{TH3as1} in the present paper
and Proposition 2 in \cite{MMas2}). This is especially important  in Theorem 1.1
of \cite{MMcol2}, where we describe the behavior after the interaction of a
solution which is as $t\to -\infty$ exactly a $2$-soliton solution.

\smallskip

Recall that the first result concerning asymptotic stability for solitons of \eqref{fkdv} was proved
by Pego and Weinstein \cite{PW}, for the power case in some weighted spaces (with exponential decay at infinity in space) under spectral assumptions, checked only 
for the nonlinearities $u^2$ and $u^3$.
This was extended by Mizumachi \cite{Mizu} under the same spectral assumptions with the condition $\int_{x>0} x^{11} u^2(x) dx<+\infty$ on the solution. 

Then, in \cite{MMjmpa} and \cite{MMarchives}, we have proved asymptotic
stability in the energy space of the solitons of \eqref{fkdv} 
in the power case respectively for $p=5$ (critical) and $p=2,$ $3$ and $4$ (subcritical).
In these papers, Theorem 2 was also the main ingredient of the asymptotic stability proof.
These results have been improved and simplified in \cite{MMnonlinearity} in the subcritical power case. The proof is direct, with no reduction to an 
Liouville property such as Theorem 2. The proof uses a  Virial identity 
which was verified only for $u^2$, $u^3$ and $u^4$ using the explicit expression of
$Q(x)$.
Finally, in \cite{yvanSIAM} the proof of the linear Liouville property (which is the main ingredient of the proof of Theorem 2) was  simplified and extended  in the power case for any $p> 1$.

\smallskip

Theorems 1 and 2 above present the first result of asymptotic  stability
of solitary waves for
\eqref{fkdv} with any nonlinearity, thus in cases where $Q_c(x)$ have no explicit
expression. In particular, the proof of Theorem 2  contains an intrinsic argument
for any solitary wave satisfying $0<c<c_*(f)$, which does not depend on a specific
potential.

We also point out that with respect to \cite{MMarchives}, the arguments to prove
Theorem 1 from Theorem 2 have been much simplified and extended. Instead of relying on
the Cauchy theory in $H^s$ for $0<s<1$ as in \cite{MMarchives}, this reduction uses
only localized energy type arguments (see proof of Proposition 4 and Appendix A).
 Moreover, the proof of Theorem 2 is direct, introducing a nonlinear dual
 problem.
 
\medskip

\noindent\emph{Remark 1.}\quad  
We focus on the case 
$Q_c>0$ (other solutions are negative and
can be addressed by changing $f(u)$ into $-f(-u)$ in equation \eqref{fkdv}).
The exponential decay assumption \eqref{decayas1} can be replaced
by an assumption of compactness of $u(t,.+ \rho(t))$ in $L^2$,
for $t\in  \mathbb{R}$ (see \cite{MMjmpa}, \cite{MMarchives}).

\medskip

\noindent\emph{Remark 2.}\quad  Note that if $Q_{c_0}$ is nonlinearly stable (in the sense that
 $\frac d{dc} {\int Q_c^2}_{|c=c_0}>0$, see Weinstein \cite{We1}), then
assumption \eqref{stableas1} can be replaced by the same assumption only at $t=0$. However,
the main point is that such a stability assumption is not needed to 
have asymptotic stability, which means that these two properties are not
related. For example, in the power case for any $p\geq 2$, $c_*(u^p)=+\infty$,
and thus Theorem 1 holds in the subcritical 
($p=2,3,4$), critical  ($p=5$) and super critical case ($p\geq 6$),
for any soliton.

In the super critical and critical  cases, the soliton is unstable
(Bona et al. \cite{BSS}, \cite{MMgafa}). In Theorem 1, we make a global assumption
on the solution (i.e. formally $u_0$ does not belong to the instable manifold of the solitons). Whether or not such solutions exist in this case is an open question, however,
the motivation of Theorem 1 in this case is to remove the possibility of any other dynamic around $Q$
(such as for example quasi-periodic solutions close to $Q$ or solutions oscillating between close solitons). In the case of the super critical nonlinear Schr\"odinger equation in dimension one,
Krieger and Schlag \cite{KRIEGERSCHLAG} have constructed a subspace of codimension 5 of initial data in which a solution close to a soliton converges to the soliton.

\medskip

\noindent\emph{Remark 3.}\quad In the case $f(u)=u^p-au^q$, where $2\leq p<q$ are integers
and $a>0$ is a constant, $c_*$ is explicit:
 $c_*=s_0^{p-1}-as_0^{q-1}$, where
$s_0=\big(\frac 1 a (\frac {q+1}{q-1})(\frac {p-1}{p+1})\big)^{\frac 1{q-p}}$.
Moreover, there is no soliton $Q_c$ for any $c>c_*$. Thus,
Theorem 1 applies to any existing  soliton in this case. For example, physical applications 
of this kind of nonlinearity in the context of the NLS equation are discussed in Sulem and Sulem \cite{SULEMSULEM}.
See also Grillakis \cite{GRILLAKIS}.

Note that the condition $f\in C^3$ can be relaxed. Indeed, all the arguments
in this paper hold for $f\in C^2$. The condition $f\in C^3$
is only assumed to obtain well-posedness of the Cauchy problem in $H^1$
by \cite{KPV89}, \cite{KPV}. More precisely, for $p\geq 3$, well-posedness 
in $H^1$ for
$f\in C^2$ follows directly from Theorem 3.6 in \cite{KPV89},
and thus Theorems 1 and 2 hold for $f\in C^2$. 
If $p=2$, one has to rely on the estimates and the norms introduced in the proof of Theorem 2.1
of \cite{KPV} for $f(u)=u^2$,
in the case $f\in C^3$ (we expect that a compactness argument should work in this case for $f\in C^2$).

\medskip

\noindent\emph{Remark 4.}\quad  It is clear that if $\frac d{dc} \int {Q_c^2}|_{c=c_0}\not = 0$ ($c$ is not critical for stability)
then $c(t)$ has a limit by Theorem 1. Our condition in Theorem 1 is more general (for example, if
$f$ is analytic, then the assumption holds unless $f(u)=u^5$). Moreover, in the case $f(u)=u^5$,
we do not expect convergence of $c(t)$ for general initial data.

\medskip

\noindent\emph{Remark 5.}\quad 
One   important tool in our analysis is a property of monotonicity of $L^2$
mass at the right in space  for solutions of the KdV equation (see Appendix A).
For the nonlinear Schr\"odinger equation, such a monotonicity property has been 
introduced in \cite{MMTnls} to prove  the stability of $N$
solitary waves for a class of suitable nonlinearities,
but so far not for proving asymptotic stability.
The question of asymptotic stability of solitary waves for the NLS equation (nonlinear 
Schr\"odinger equation) is mostly open, see results for special nonlinearities
by Buslaev and Perleman \cite{BP}, Perelman \cite{P} and Rodnianski, Schlag and Soffer \cite{RSS}. It is a promising direction of research.

\medskip

\noindent\emph{Remark 6.}\quad 
In the integrable case ($f(u)=u^2$), using the inverse scattering transform,
a general decomposition result has been proved by  Eckauss and Schuur \cite{ES}:
any smooth ($C^4$) and sufficiently decaying solution decomposes as $t\rightarrow +\infty$ in a finite sum of solitons plus a dispersive part that converges to zero in some sense. This implies the result of Theorem 1 for this nonlinearity and such initial data.
Such questions for the integrable NLS equation (cubic NLS equation in one space dimension) are open.

\medskip

The paper is organized as follows. In Section 2, we recall some prelimary results concerning
solutions of \eqref{ellipticas1}. In Section 3, we prove Theorem 2 and in Section 4, we prove
Theorem 1. Section 5 is devoted to the multi-soliton case.

\medskip

\noindent\textbf{Acknowledgements.} We wish to thank the referees for their 
useful comments.

\section{Preliminary results}

\subsection{Stationary equation \eqref{ellipticas1}}
First, we recall the necessary and sufficient condition for existence
of a solution of \eqref{ellipticas1}, and some properties of the solution.
Let $f$ be $C^2$ and satisfy \eqref{surfas1} (so that for any $c>0$,
$\frac c2 s^2-F(s)>0$ for small positive $s$).

\begin{claim}\label{BLas1} Let $c>0$.
There exists a nontrivial solution $Q_c\in H^1(\mathbb{R})$ $(Q_c(x_0)>0$ for $x_0\in \mathbb{R})$
of \eqref{ellipticas1}
if and only if there exists $s_0>0$ the smallest positive zero of $s\mapsto \frac c2 s^2-F(s)$ and $s_0$
satisfies $cs_0-f(s_0)<0$. 

Moreover, $Q_c$ is $C^4$, unique up to translation and can be chosen so that 
$Q_c(0)=s_0$,
$Q_c(x)=Q_c(-x)$, $Q_c'(x)>0$ for all $x>0$. 
Finally, there exists $K>0$ such that
\begin{equation}\label{decayQas1}
\forall x\in \mathbb{R},\quad
\frac 1K e^{-\sqrt{c} |x|} \leq Q_c(x)\leq K e^{-\sqrt{c} |x|},
\quad  |Q_c'(x)|\leq K e^{-\sqrt{c} |x|}.
\end{equation}
\end{claim}

\noindent\emph{Proof.} We refer to  Berestycki and Lions \cite{BL}, Theorem 5 and Remark 6.3 in section 6
for the proof of these results.

\medskip

By assumption \eqref{surfas1} and Claim \ref{BLas1}, there exists $\bar c>0$ such that for any
 $0<c<\bar c$, $Q_c$ exists with $\|Q_c\|_{L^\infty}\to 0$ as $c\to 0$.
 Thus we may define
$$c_*=\sup\{c>0 \text{ such that }
\forall c'\in ]0,c[, \text{ $\exists$ $Q_{c'}$   positive solution of 
\eqref{ellipticas1}}\}.$$
In the power case, we have $c_*=+\infty$ by scaling property.
Note also that if $c_*<+\infty$ then from Claim \ref{BLas1}
there exists no nontrivial solution of \eqref{ellipticas1} for $c={c_*}$.

Let us give another characterization of $c_*$, which is the one used in the
proofs.

\begin{claim}\label{equivbis}
A unique even positive solution $Q_c$ of \eqref{ellipticas1} exists and satisfies
\begin{equation}\label{cccbisas1}
\forall x\in \mathbb{R},\quad
Q_c(x) f(Q_c(x)) - 2 F(Q_c(x))>0
\end{equation}
if and only if $0<c<c_*$.
\end{claim}

Note that this property
is related to the Palais-Smale condition for the corresponding variational
problem in dimension greater or equal than $2$.

\medskip

\noindent\emph{Proof.} 
First, let $c>0$ and assume the existence of $Q_c>0$ solution of \eqref{ellipticas1} satisfying
\eqref{cccbisas1}.  Let $s_c=Q_c(0)$. Since $Q_c(\mathbb{R})=(0,s_c]$,
by \eqref{cccbisas1}, we have:
\begin{equation}
\label{wqwq}
\forall s\in (0,s_c],\quad sf(s)-2F(s)>0.
\end{equation}
Let $0<c'<c$. 
Let us prove that there exists a positive solution of \eqref{ellipticas1} for
$c'$.
Since
$\frac {c'}2 s_c^2 -F(s_c) < \frac {c}2 s_c^2 -F(s_c)=0$
and \eqref{surfas1}, there exists
$0<s_{c'}<s_c$ the first zero of $\frac {c'} 2 s^2 -F(s)$, and by 
\eqref{wqwq},   $s_{c'} f(s_{c'})-2 F(s_{c'})>0$.
Together with $\frac {c'}2 s_{c'}^2 -F(s_{c'})=0$ this implies that 
$c' s_{c'} -f(s_{c'})<0$ and thus by Claim \ref{BLas1},
there exists $Q_{c'}$ solution of \eqref{ellipticas1} with $c=c'$. Since
$0<c'<c$ is arbitrary, we have proved $0<c<c_*$.

Second, let us consider $0<c<c_*$. For the sake of contradiction
assume that for some $0<s\leq Q_c(0)$, 
$sf(s)-2F(s)\leq 0$. Let $0<s_1\leq Q_c(0)$ be the smallest such $s$,
so that by \eqref{surfas1},
$s_1 f(s_1)-2 F(s_1)=0$ and $sf(s)-2F(s)>0$ on $(0,s_1)$.
Let $c'=\frac {F(s_0)}{\frac 12 s_0^2}$. Since $s\mapsto \frac {F(s)}{\frac 12 s^2}$
is strictly increasing on $[0,s_1]$, $s_1$ is the first zero of
$s\mapsto \frac {c'} 2  s^2-F(s)$. Using  $s_1 f(s_1)-2 F(s_1)=0$, we obtain that 
$c's_1-f(s_1)=0$, which implies that equation \eqref{ellipticas1} has no solution
for $c=c'$, a contradiction with the definition of $c_*$.

\subsection{Linearized operator around $Q_c$}
 Let $\varphi(x)$ be a $C^2$ even function such that
$0\leq \varphi\leq 1$, $|\varphi'|+|\varphi''|\leq 4$ on $\mathbb{R}$,
$\varphi\equiv 1$ on $[0,1]$, $\varphi\equiv 0$ on $[2,+\infty)$. 
Let $\langle f,g\rangle$ denote the $L^2$ scalar product of $f$ and $g$.
We consider the linearized operator around $Q_{c_0}$:
\begin{equation}
\label{linearized}
\mathcal{L}_{c_{0}}= -\partial_x^2 + {c_{0}} - f'(Q_{c_{0}}).
\end{equation}

\begin{claim}\label{CL1as1} Let $0<{c_{0}}<c_*(f)$. The following properties hold
\begin{enumerate}
  \item There exist unique $\lambda_0>0$, $\tilde \chi_{c_{0}}\in H^1(\mathbb{R})$,
  $\tilde \chi_{c_{0}}>0$ such that
  $\mathcal{L}_{c_{0}} \tilde\chi_{c_{0}}=-\lambda_0 \tilde\chi_{c_{0}}$,
  $\langle \tilde \chi_{c_0},\tilde \chi_{c_0}\rangle=1$.
  \item For all $u\in H^1$, $\mathcal{L}_{c_{0}} u =0$ is equivalent to $u=\lambda Q_{c_{0}}'$,
  $\lambda\in \mathbb{R}$.
    \item For all $h\in L^2$, if $\langle h,Q'_{c_0}\rangle=0$ then there exists a 
    unique $u\in H^2$ such that $\langle u,Q'_c\rangle=0$ and $\mathcal{L}_c u=h$.
  
  Moreover,  $\mathcal{L}_{c_0} S_{c_0}=-Q_{c_0}$ where $S_{c_0}=\frac d{dc} {Q_{c}}_{| c={c_0}}$.
     \item There exist $B,\lambda_1, \sigma_1>0$ such that 
  for all $c \in [c_0-\sigma_1,c_0+\sigma_1]$,
  the function $  \chi_c(x)=\tilde\chi_c(x)\varphi(\frac x B)$
  satisfies
  \begin{equation}\label{avantPOS}
  \int \chi_c Q_c>0,\quad 
  \frac {\lambda_0}2\leq 
  -\frac {\langle\mathcal{L}_c \chi_c,\chi_c\rangle}{\langle \chi_c,\chi_c\rangle}
  \leq {\lambda_0},  \end{equation}
  \begin{equation}\label{POS}
  \forall u\in H^1(\mathbb{R}),
  \quad\int uQ_c'=\int u  \mathcal{L}_c\chi_c=0\quad \Rightarrow \quad
  \langle\mathcal{L}_c u,u\rangle\geq \lambda_1 \langle u,u\rangle.
  \end{equation}
\end{enumerate}
\end{claim}

\noindent\emph{Proof.} The first three properties follow from classical arguments.
See Weinstein \cite{We2}, proof of Proposition~2.8b for $N=1$ and proof of Proposition 2.10.
Note that letting $S_{c_0}=\frac d{dc} {Q_{c}}_{| c={c_0}}$, then
by differentiating the equation of $Q_{c}$ with respect to $c$, we obtain $\mathcal{L}_{c_0} S_{c_0}=-Q_{c_0}$.
Note also that  $\tilde \chi_c(x)\leq K e^{-\sqrt{c} |x|}$ and 
\eqref{POS} holds for $\tilde \chi$.

Now, let 
$\chi_c(x)=\tilde\chi_c(x)\varphi(\frac x B)$. By index theory of quadratic form, it is enough to check \eqref{avantPOS}.
We have  $\chi_c\geq 0$ and $\chi_c\not \equiv 0$, so that $\int \chi_c Q_c>0$,  
$\langle \chi_c,\chi_c\rangle = 1 + O(e^{-\sqrt{c}B})$ and 
$\mathcal{L}_c \chi_c=(\mathcal{L}_c \tilde \chi_c)
\varphi(\frac x B) -\frac 2B \varphi'(\frac x B) \tilde \chi_c'(\frac x B) - \frac {1}{B^2} \varphi''(\frac x B)
\tilde \chi_c$ so that
$\langle \mathcal{L}_c \chi_c,\chi_c\rangle = -\lambda_0  + O(e^{-\sqrt{c} B})$.
Thus, \eqref{avantPOS} follows by taking $B$ large enough.

\medskip

From now on, $B$ is fixed to such value. Note that $\chi_c$ has support in $[-2B,2
B]$, uniform in $c \in [c_0-\sigma_1,c_0+\sigma_1]$.

\subsection{Decomposition of a solution close to $Q_{c_0}$}

\begin{lemma}[Modulation of a solution close to $Q_{c_0}$]\label{DECOMPas1}
Let $0<c_0<c_*$.
There exist $K_0>0$ and $\alpha_0>0$ such that for any $0<\alpha<\alpha_0$
and $T_0>0$, if $u(t)$ solution of \eqref{fkdv} satisfies
\begin{equation}
\forall t\in [0,T_0],\quad \inf_{r\in \mathbb{R}}\|u(t)-Q_{c_0}(.-r)\|_{H^1}
\leq \alpha,
\end{equation}
then there exist $c(t)>0$, $\rho(t)\in C^1([0,T_0])$ such that
\begin{equation}
\label{decomptt}
\eta(t,x)=u(t,x)-Q_{c(t)}(x-\rho(t)),
\end{equation}
satisfies, for all $t\in [0,T_0],$
\begin{align}
\label{propeps}
    &   \int \tilde \chi_{c(t)}(x-\rho(t)) \eta(t,x)dx=
    \int Q'_{c(t)}(x-\rho(t)) \eta(t,x)dx=0, \\
    &	|c(t)-c_0|+\|\eta(t)\|_{H^1}\leq K_0 \alpha,
    \\
    &  |c'(t)|+|\rho'(t)-c(t)|\leq K_0 \left(\int \eta^2(t,x) e^{-|x-\rho(t)|}dx\right)^{\frac 12}. \label{decompttbis}
\end{align}
\end{lemma}

\noindent\emph{Proof.} This is a standard application of the implicit function theorem.
See for example \cite{MMgafa}, Proposition 1 for details.
Note that 
  $\frac d{dc'} {Q_{c'}}_{|c'=c}=-S_c$ 
and $\frac d{dx'} {Q_c(x+x')}_{|x'=0}=Q_c'(x)$. Thus,
the nondegeneracy conditions  are (by Claim \ref{CL1as1}),
$$\int S_c   \tilde \chi_c= - \frac 1{\lambda_0} \int  \mathcal{L}_c (S_c) {\tilde \chi_c}
=\frac 1{\lambda_0} \int   Q_c {\tilde \chi_c}>0,\quad
\int \tilde \chi_c Q_c'=0,$$   
$$\int Q_c' S_c=0,\quad \int (Q_c')^2>0.$$

\section{Rigidity results}

This section is devoted to the proof of the rigidity theorem (Theorem 2),
see Section 3.2.
First, in section 3.1, we give a linear version of the result
to present the main idea in the simplest case.

\subsection{Linear Liouville property}

In this section, under the assumptions of Theorem 2,
we prove a rigidity result for $H^1$ solutions of
the following linearized equation
\begin{equation}\label{eqeta}
\partial_t \eta  = \partial_x (\mathcal{L}_{c_0} \eta), \quad (t,x)\in \mathbb{R}\times \mathbb{R},\quad
\text{where} \quad \mathcal{L}_{c_0} \eta= -\eta_{xx}+ {c_0} \,  \eta -f'(Q_{{c_0}}) \eta.
\end{equation}
Note that the arguments of Lemma 9 in \cite{MMjmpa} (based on linear estimates of Kenig, Ponce and Vega \cite{KPV}) 
prove that the Cauchy problem for \eqref{eqeta} is globally well-posed in $H^1(\mathbb{R})$
(in a certain sense). By $H^1$ solution, we mean a solution constructed in this way. Any such solution can be approached by regular solutions
 which allows to justify formal computations.

\begin{proposition}[Linear Liouville property]\label{PROPas1}
Let $0<{c_0}<c_*(f)$. 
Let  $\eta \in C(\mathbb{R},H^1(\mathbb{R}))$ be solution of \eqref{eqeta}.
Assume that there exist $K>0$, $\sigma>0$ such that
\begin{equation}\label{dsds}
\forall (t,x) \in \mathbb{R}\times \mathbb{R},\quad
|\eta(t,x)|\leq K e^{-\sigma |x|}.
\end{equation}
Then, there exists $b_0\in \mathbb{R}$ such that for all $t\in \mathbb{R}$,
$\eta(t)\equiv b_0 Q_{{c_0}}'$.
\end{proposition}

\noindent\emph{Remark.} Note that since $Q_{{c_0}}'$ verifies $\mathcal{L}_{c_0} Q_{c_0}'=0$
and has exponential decay, $\eta(t)\equiv b_0 Q_{{c_0}}'$ is solution of \eqref{eqeta}--\eqref{dsds}.

\medskip

Let $\eta(t)$ be an $H^1$ solution of  \eqref{eqeta}
satisfying \eqref{dsds}. As in \cite{yvanSIAM}, we introduce a dual problem related to $\eta$.
\begin{lemma}[Introduction of the dual problem]\label{SURV}
Let 
$$v(t,x)=\mathcal{L}_{c_0} \eta(t,x) + \alpha(t) Q_{c_0}\quad \text{where}\quad
\alpha(t)=- \frac {\int \eta \mathcal{L}_{c_0} \chi_{c_0}}{\int \chi_{c_0}Q_{c_0}}.$$
Then, $v\in C(\mathbb{R}, H^1(\mathbb{R}))$ and
$v(t)$ satisfies:
\begin{enumerate}
\item Equation of $v$.
\begin{equation}\label{eqvv}
  \partial_t   v = \mathcal{L}_{c_0} \left( \partial_x   v \right) +\alpha'(t) Q_{{c_0}},
  \quad (t,x)\in \mathbb{R}\times \mathbb{R}.
 \end{equation}
\item Exponential decay. There exists $K>0$ such that
\begin{equation}\label{expov}	
 \forall (t,x)\in \mathbb{R}\times\mathbb{R},\quad 
 \left|v(t,x) \right|
 \leq K e^{-\frac {\sqrt{c_0}} 8 |x|}. \end{equation}
\item Orthogonality relations.
\begin{equation}\label{orthov}
\forall t\in \mathbb{R},\quad
\int v(t,x) \chi_{{c_0}}(x) dx=\int v(t,x) Q_{{c_0}}'(x) dx=0.
\end{equation}
\item Virial identity on the dual problem. Let
\begin{equation}\label{vireq1}
\mu_{c_0}(x)= -\frac {Q_{{c_0}}'(x)}{Q_{c_0}(x)}\quad \text{then}\quad
\frac 12 \frac d{dt} \int v^2(t,x) \mu_{c_0}(x) dx
= \int  \partial_x v \, \mathcal{L}_{c_0} \left( v \mu_{c_0} \right).
\end{equation}
\end{enumerate}
\end{lemma}
 
\noindent\emph{Proof of Lemma \ref{SURV}.}
1. We have $\eta_t = \partial_x (\mathcal{L}_{c_0} \eta) = v_x -\alpha(t) Q_{c_0}'$, thus
using $\mathcal{L}_{c_0}Q_{c_0}'=0$, we obtain
$$v_t = \mathcal{L}_{c_0} \eta_t  + \alpha'(t) Q_{c_0} 
=\mathcal{L}_{c_0} v_x  + \alpha'(t) Q_{c_0}.$$

\smallskip

2. From  monotonicity arguments on $\eta(t)$ and on $v(t)$, we claim
that there exists $K>0$ such that for all $t\in \mathbb{R}$,
\begin{equation}\label{monov}
\int (v_x^2+v^2)(t) \exp\left(\frac {\sqrt{c_0}} 4 |x|\right) dx \leq K.
\end{equation}
(proof in Appendix A).
By $\|v \exp\left(\frac {\sqrt{c_0}} 8 |x|\right)\|_{L^\infty}
\leq K \|v \exp\left(\frac {\sqrt{c_0}} 8 |x|\right)\|_{H^1}$,
this implies \eqref{expov}. The fact that
$v\in C(\mathbb{R}, H^1(\mathbb{R}))$ then follows from the equation.

\smallskip

3. By the choice of $\alpha(t)$ and $\mathcal{L}_{c_0} Q_{c_0}'=0$, we have
$$\int v \chi_{c_0}=\int \mathcal{L}_{c_0} \eta \, \chi_{c_0} + \alpha(t)Ê\int Q_{c_0}\chi_{c_0}=0,
\quad \int v Q'_{c_0}=\int \mathcal{L}_{c_0} \eta \, Q_{c_0}' + \alpha(t)Ê\int Q_{c_0}Q_{c_0}'=0.$$

\smallskip

4. From the equation of
$v$ and $\int Q_{{c_0}} v(t) \mu_{c_0}= -\int  v(t) Q_{{c_0}}' =0$,
 we have
\begin{equation*}
\frac 12 \frac d{dt} \int v^2(t) \mu_{c_0} =   \int v\partial_t v \,  \mu_{c_0}
 =   \int \mathcal{L}_{c_0}( \partial_x v)  v \, \mu_{c_0}
  +  \, \alpha'(t) \int Q_{{c_0}} v \, \mu_{c_0}
=  \int  \partial_x v \, \mathcal{L}_{c_0}( v\,  \mu_{c_0}).
\end{equation*}

\medskip

Now, we claim the following.

\begin{lemma}[Positivity of the quadratic form]\label{VIRLIN}
For any $0<c<c_*(f)$, there exists $\lambda_2(c)>0$ continuous such that
\begin{equation}\label{surg100}
\forall x \in \mathbb{R},\quad
\frac {\lambda_2}{\cosh^{p-1}(\sqrt{c} x)} \leq \mu_c'(x) \leq
\frac 1{\lambda_2}\frac 1{\cosh^{p-1}(\sqrt{c} x)},
\end{equation}
\begin{equation}\label{virlin1}
\forall  w \in H^1,\quad 
 -\int  \partial_x w \, \mathcal{L}_{c} \left( w \mu_c \right) 
 =\frac 32 \int (\partial_x (\tfrac {w}{Q_c}))^2 Q_c^2 \mu_c'\geq \lambda_2 \int w^2 \mu_c' -\frac 1{\lambda_2} \left(\int w \chi_c\right)^2.
\end{equation}
\end{lemma}

\noindent\emph{Proof of Lemma \ref{VIRLIN}.} 
First, by $Q''_{c}={c} Q_{c} - f(Q_{c})$ and $(Q_{c}')^2 = {c} Q_{c}^2 - 2F(Q_{c})$,
we have by Claim \ref{equivbis} and  $0<{c}<c_*$, $x\neq 0$,
\begin{equation}\label{pourthetap}
\mu_c'=\frac 1{Q_{c}^2} ((Q_{c}')^2 - Q_{c} Q_{c}'')=\frac 1{Q_{c}^2} (Q_{c} f(Q_{c}) -2 F(Q_{c}))>0
\end{equation}
and we obtain \eqref{surg100} from \eqref{surfas1} and
\eqref{decayQas1} and continuity arguments.
%\frac {\delta}{\cosh(\sqrt{c} x)} \leq Q_c(x) \leq
%\frac 1{\delta}\frac 1{\cosh(\sqrt{c} x)},
%\quad

Next, let $z= \frac w{Q_{c}}$ so that $w\mu_c=-zQ_{c}'$. 
We claim 
\begin{equation}
\label{encorevir}
-\int  \partial_x w \, \mathcal{L}_{c} \left( w \mu_c \right) 
=\frac 32\int (\partial_x z)^2 Q_c^2 \mu_c'.
\end{equation}
Using 
$$\mathcal{L}_{c} \left( z Q_{{c}}' \right)=
z \mathcal{L}_{c} Q_{c}' - 2 \partial_x z \,Q_{{c}}'' - \partial_x^2 z \,Q_{{c}}'=- 2 \partial_x z \,Q_{{c}}'' - \partial_x^2 z \,Q_{{c}}',$$
we have
\begin{equation*}\begin{split} &
-\int  \partial_x w \, \mathcal{L}_{c} \left( w \mu_c \right)   = \int \partial_x (Q_{{c}}z) \mathcal{L}_{c} \left( z Q_{{c}}' \right)
=\int (Q_{{c}}' z + Q_{{c}} \partial_x z) (- 2 \partial_x z \,Q_{{c}}'' - \partial_x^2 z \,Q_{{c}}')\\
& = \int \left(z (Q_{c}'Q_{c}'')' + (\partial_x z)^2 (Q_{c}')^2
-\tfrac 12 z^2 ((Q_{c}')^2)''- 2 (\partial_x z)^2 Q_{c} Q_{c}'' 
+\tfrac 12 (\partial_x z)^2 (Q_{c} Q_{c}')'\right)
\\ &=\frac 32 \int (\partial_x z)^2 ((Q_{{c}}')^2 - Q_{{c}} Q_{{c}}'')
=\frac 32 \int (\partial_x z)^2 Q_c^2 \mu_c',
\end{split}
\end{equation*}
by \eqref{pourthetap}, which proves \eqref{encorevir}.

\medskip

Let $Z(x)=   z(x)\,  {\cosh^{-\frac {p+1} 2} (\sqrt{c} x)}$. By \eqref{surg100} and direct computations,
 we have ($\delta>0$)
\begin{equation}\label{Ltilde}
 \int (\partial_x z)^2 \mu_c' Q_c^2
\geq\delta \int (\partial_x z)^2 {\cosh^{-p-1}(\sqrt{c} x)} 
= \delta\langle \tilde{\mathcal{L}}_c Z, Z\rangle,
\end{equation}
where  $\tilde{\mathcal{L}}_c Z = -Z_{xx}+\frac c4 \left(  {p+1}  \right)^2 Z
- \frac c4 {(p+1)(p+3)}  Z  \cosh^{-2}(\sqrt{c} x).$
By \eqref{Ltilde}, $\tilde{\mathcal{L}}_c$ is a nonnegative operator, 
with first eigenvalue $0$ associated to the function
 ${\cosh^{-\frac {p+1} 2} (\sqrt{c} x)}$.
From standard arguments, since the function $Q_c \chi_c  {\cosh^{\frac {p+1} 2} (\sqrt{c} x)}\geq 0$
is nonnegative, not zero
and belongs to $L^2$ (this is where we use that $\chi_c$ is compactly supported),
there exists $\lambda>0$ such that 
$$-\frac 23\int \partial_x  w\mathcal{L}_c (w \mu_c)=
\int (\partial_x z)^2 \mu_c' Q_c^2
\geq \delta\langle \tilde{\mathcal{L}}_c Z, Z\rangle
 \geq \lambda \int Z ^2 
- \frac 1{\lambda} \left(\int Z\, Q_c \chi_c  {\cosh^{\frac {p+1} 2} (\sqrt{c} x)}\right)^2.
$$
Since $w=zQ_c= Z  Q_c\,  {\cosh^{\frac {p+1} 2} (\sqrt{c} x)} $,  from  \eqref{surg100}, we obtain ($\lambda_2>0$)
$$-\int \partial_x  w\mathcal{L}_c (w \mu_c) \geq \lambda_2 \int w^2 \mu_c' 
- \frac 1{\lambda_2} \left(\int w \,   \chi_c \right)^2.
$$

\medskip

\noindent\emph{Proof of Proposition \ref{PROPas1}.}
By \eqref{vireq1}, Lemma \ref{VIRLIN},  and \eqref{orthov}, we have
\begin{equation}\label{virfinal}
-\frac 12\frac d{dt}Ê\int v^2(t) \mu_{c_0} \geq \lambda_2 \int v^2 \mu_{c_0}'.
\end{equation}
Since $|\mu_{c_0}(x)|\leq C$ on $\mathbb{R}$ and $v(t)$ is uniformly bounded in time in $L^2$,
  $\lim_{t\to \pm \infty} \int v^2(t) \mu_{c_0}=l_{\pm}$ exist
 and  by integrating \eqref{virfinal},
\begin{equation}\label{intee}
\int_{-\infty}^{+\infty} \int v^2(t,x) \mu'_{c_0}(x) dx dt
\leq \frac 1 {2\lambda_2} (l_- -l_+)
<+\infty.
\end{equation}
By \eqref{surg100}, it follows that for a sequence $t_n\to +\infty$, we have 
 $v(t_n)\to 0$ in $L^2_{loc}(\mathbb{R})$ 
 and thus by \eqref{expov}, 
$v(t_n)\to 0$  in $L^2(\mathbb{R})$ as $n\to +\infty$ and 
 $l_+=0.$
Similarly, 
$l_-=0$.
Thus, by \eqref{intee} and $v\in C(\mathbb{R},H^1)$, we obtain
$$\forall (t,x)\in \mathbb{R}\times\mathbb{R},\quad v(t,x)= 0.$$

It follows that $\mathcal{L}_{c_0} \eta(t)= - \alpha(t)  Q_{c_0}$.
Thus, by Claim \ref{CL1as1}, we obtain, for some bounded
function $\beta(t)$,
$$\eta(t)=\alpha(t) S_{c_0} + \beta(t) Q_{c_0}'.$$
By the equation of $\eta(t)$ \eqref{eqeta},
and the orthogonality of $S_{c_0}$ and $Q_{c_0}'$, we obtain
$\beta'(t)=-\alpha(t)$ and $\alpha'(t)=0$. Since $\beta(t)$ and $\alpha(t)$ are bounded, we deduce
$\alpha(t)\equiv 0$ and $\beta(t)\equiv b_0$.

\subsection{Nonlinear Liouville property - Proof of Theorem 2}

The proof of Theorem \ref{PROP2as1} follows the same lines as
the proof Proposition \eqref{PROPas1}.
Consider now $u(t)$ as in Theorem 2.
We first decompose $u(t,x)$ similarly as in  Lemma \ref{DECOMPas1}, using modulation theory.
We obtain, for all $t\geq 0$,
\begin{equation}\label{decomp}
\eta(t,x)=u(t,x+\rho(t))-Q_{c(t)}(x),
\end{equation}
where $c(t)$, $\rho(t)$ are $C^1$ functions chosen so that 
\begin{equation}\label{orthoeta}
\int \eta(t,x)\mathcal{L}_{c(t)}\chi_{c(t)}(x) dx=\int \eta(t,x) Q_{c(t)}'(x) dx=0.
\end{equation}
(The nondegeneracy conditions in this case
are $\int S_c \mathcal{L}_c \chi_c= \int  \mathcal{L}_c (S_c) {\chi_c}
=- \int   Q_c {\chi_c}<0$ and $\int (Q_c')^2>0$.)
Recall that 
\begin{equation}\label{etapetit}
\|\eta(t)\|_{H^1}+|c(t)-c_0|\leq K \alpha_0.
\end{equation} 
Thus, we can choose $\alpha_0>0$ small enough so that,
for all $t\geq 0$, $c(t)\in [c_0-\sigma_0,c_0+\sigma_0]\subset (0,c_*)$,
for $\sigma_0>0$ small enough so that Claim \ref{CL1as1} and Lemma \ref{VIRLIN} apply to $c=c(t)$.

% We recall $S_c=\frac d {d\tilde c } {Q_{\tilde c}}_{|\tilde c=c}$ and $\mathcal{L}_c S_c=-Q_c$.

As for the linear equation, we introduce a dual problem.
 
\begin{lemma}[Dual problem for the nonlinear equation]\label{CL100asz1}
Let 
\begin{equation*}\begin{split}
v(t,x)&=-\eta_{xx} + c\eta - (f(Q_c+\eta)-f(Q_c))\\
&=\mathcal{L}_c  \eta - (f(Q_c+\eta)-f(Q_c)-f'(Q_c)\eta).
\end{split}\end{equation*}
Then, $v\in C(\mathbb{R},H^1(\mathbb{R}))$ and $v(t)$ satisfies
\begin{enumerate}
  \item Equation of $v$.
  \begin{equation}\label{eqvas1}\begin{split}
v_t&= -v_{xxx}+c v_x - v_x f'(Q_c+\eta) + (\rho'-c) v_x + c' (Q_c+\eta)\\
&= \mathcal{L}_c (v_{x}) - v_x (f'(Q_c+\eta)-f'(Q_c)) + (\rho'-c) v_x + c' (Q_c+\eta).\end{split}\end{equation}
	\item Exponential decay. There exists
  $K>0$ such that,
	\begin{equation}
\label{decayvNL}
\forall (t,x)\in \mathbb{R}\times \mathbb{R},\quad
|\eta(t,x)|+|v(t,x)|\leq K e^{-\frac {\sqrt{c_0}}8 |x|}.
\end{equation}
  \item Estimates and almost orthogonality relations. There exists
  $K>0$ such that, $\forall t\in \mathbb{R}$,
\begin{equation}
\label{orthovNL}
|c'|  + |\rho'-c|\leq K \|\eta\|_{L^2},\quad  \left| \int v Q_c'\right| +
\left| \int v \chi_c\right|\leq K \|\eta\|^2 _{L^2},\quad 
\|\eta\|_{L^2}       \leq K \| v \|_{L^2}.
\end{equation}
\item Virial type estimates.
There exists $\lambda_3, B>0$ such that, $\forall t\in \mathbb{R}$,
\begin{align}
    & -\frac 12 \frac d{dt} \int v^2  \mu_c \geq \lambda_3 \int v^2 \mu_c' - 
    {\frac 1 {\lambda_3}} \|v\|_{H^1}^2 \|\eta\|_{L\sp2}
    			\label{VIR1as1}\\
    & -\frac 12 \frac d{dt} \int x v^2 \geq \lambda_3 \int (v_x\sp2 +v\sp2)   - {\frac 1 {\lambda_3}} \int_{|x|\leq B} v^2  
    			\label{VIR2as1}
			\end{align}
\end{enumerate}
\end{lemma}

\noindent\emph{Remark.} Note that at the first order, we have $v(t)\sim \mathcal{L}_{c}\eta(t)$,
$\int v \chi_c\sim 0$ and $\int v Q_c'\sim 0$ as in the proof of the linear Liouville property.

\medskip

\noindent\emph{Proof of Lemma \ref{CL100asz1}.}
1. First, we write the equation of  $\eta(t)$, from \eqref{decomp}, \eqref{fkdv}
and \eqref{ellipticas1}
\begin{equation}\label{eqetaNL}\begin{split}
 \eta_t
 & = u_t +\rho' u_x - c' S_c = -(u_{xx}+f(u))_x + \rho' u_x - c'S_c
 \\ & = (-\eta_{xx} + c\eta - (f(Q_c+\eta)-f(Q_c)))_x + (\rho'-c) (Q_c+\eta)_x -c' S_c
 \\ & = v_x + (\rho'-c) (Q_c+\eta)_x -c' S_c,
\end{split} \end{equation}
where $v=-\eta_{xx} + c\eta - (f(Q_c+\eta)-f(Q_c)).$
Now, we compute $v_t$:
\begin{equation*}\begin{split}
v_t 
& = -\eta_{txx}+c\eta_t -\eta_t f'(Q_c+\eta) + c'\eta - c' S_c (f'(Q_c+\eta)-f'(Q_c))
\\ & = - v_{xxx} + cv_x - v_x f'(Q_c{+}\eta)
+ (\rho'-c) (-(Q_c {+} \eta)_{xxx} + c (Q_c{+}\eta)_x - (Q_c{+} \eta)_x f'(Q_c{+}\eta))\\ &
- c' (-S_{cxx} + cS_c -S_c f'(Q_c{+}\eta)) + c'\eta- c' S_c (f'(Q_c{+}\eta)-f'(Q_c)).
\end{split}\end{equation*}
Since 
$v_x= -\eta_{xxx}+c \eta_x - \eta_x f'(Q_c+\eta) - Q_c'(f'(Q_c+\eta)-f'(Q_c))$,
we obtain
\begin{equation*} 
v_t 
 = - v_{xxx} + cv_x - v_x f'(Q_c{+}\eta)
+ (\rho'-c) (v_x+ \mathcal{L}_c Q_c') - c' \mathcal{L}_c S_c + c'\eta.
 \end{equation*}
Thus, by
$\mathcal{L}_c Q_c'=0$ and $\mathcal{L}_c {S}_c=-Q_c$ (see Claim \ref{CL1as1}), we obtain \eqref{eqvas1}.

\smallskip

2. By monotonicity arguments, we claim that 
there exists $K>0$ (independent of $\alpha_0$)
such that for all $t\in \mathbb{R}$,
\begin{equation}\label{monovNL}
\int (v_x^2+v^2)(t) \exp\left(\frac {\sqrt{c_0}} 4 |x|\right) dx \leq K.
\end{equation} 
See the proof of   \eqref{monovNL} in Appendix A. Note that \eqref{monovNL}
 implies \eqref{decayvNL} (see proof of Lemma \ref{SURV}).

\smallskip

3. By classical arguments (multiply \eqref{eqetaNL}
 by $\chi_c$ (respectively, by $Q_c'$)
and integrate on $\mathbb{R}$), we obtain $|c'| + |\rho'-c|  \leq K \| \eta \|_{L^2}$.
See \cite{MMgafa} for example.

Next,
$\int v Q_c'=\int \mathcal{L}_c \eta Q_c' - \int (f(Q_c+\eta)-f(Q_c)-f'(Q_c)\eta) Q_c'$
and since $\mathcal{L}_c Q_c'=0$ and $|f(Q_c+\eta)-f(Q_c)-f'(Q_c)\eta|\leq K \eta^2$
($f$ is $C^2$), we obtain
$\left| \int v Q_c'\right| \leq K \int \eta^2$.
Since  
$\int \eta \mathcal{L}_c\chi_c=0$, 
$\int v \chi_c = \int (\mathcal{L}_c \eta) \chi_c- \int (f(Q_c+\eta)-f(Q_c)-f'(Q_c)\eta) \chi_c
$ implies $\left| \int v \chi_c\right| \leq K \int \eta^2$.

By Claim \ref{CL1as1} and \eqref{orthoeta}, we have $\langle \mathcal{L}_c \eta,\eta\rangle \geq \lambda_1\int \eta^2$.
Thus, since $f$ is $C^2$,
$$\langle v,\eta\rangle = \langle \mathcal{L}_c \eta,\eta\rangle
- \int (f(Q_c+\eta){-}f(Q_c){-}f'(Q_c)\eta) \eta
\geq \lambda_1 \int \eta^2 - K \|\eta\|_{L^\infty} \int \eta^2
\geq \tfrac 12 \lambda_1 \int \eta^2,$$
for $\alpha_0$ small enough using \eqref{etapetit}. Thus 
$\int \eta^2 \leq K \int v^2$ by Cauchy-Schwartz inequality.

\medskip

4. \emph{Proof of \eqref{VIR1as1}.}
By the equation of $v$
\begin{equation*}
\begin{split} &
-\frac 12 \frac d{dt} \int   v^2 \mu_{c} = -\int v_t v \mu_{c} - \frac {c'}2  \int v^2  \frac {d\mu_{c}}{dc} = -\int v_x \mathcal{L}_c (v\mu_{c}) +R_1
\quad\text{where}
\\ & R_1
=  \int (f'(Q_c+\eta)-f'(Q_c)) \mu_{c} v_x v 
+ \frac 12 (\rho'-c)\int  v^2 \mu_{c}'  + c' \int v Q_c' 
-c' \int \eta v \mu_c
-  \frac {c'}2 \int v^2  \frac {d\mu_{c}}{dc}.
\end{split}
\end{equation*}
By Lemma \ref{VIRLIN} and \eqref{orthovNL}, we have
$$
-\int v_x \mathcal{L}_c (v\mu_{c}) \geq \lambda_2 \int v^2 \mu'_c - K \|\eta\|_{L^2}^4
\geq \lambda_2 \int v^2 \mu'_c - K \|\eta\|_{L^2}\|v\|_{L^2}^2.
$$
Now, we prove $|R_1|\leq K \|\eta\|_{L^2} \|v\|_{H 1}^2$ and \eqref{VIR1as1}
will follow.

Since  $f$ is $C^2$, we have  $|f'(Q_c+\eta)-f'(Q_c)|\leq K |\eta|$
and so
$$\left|\int (f'(Q_c+\eta)-f'(Q_c)) \mu_{c} v_x v \right|
\leq K \|v\|_{L^\infty} \int |\eta| |v_x|\leq K \|v\|_{H^1}^2 \|\eta\|_{L^2}.$$
By \eqref{orthovNL} and since
$\mu_c$, $\mu'_c$, $\frac {d\mu_{c}}{dc}$ are bounded, we have
 $$|(\rho'-c)\int  v^2 \mu_{c}'|+|c' \int v Q_c'|+|c' \int \eta v \mu_c|+|\frac {c'}2 \int v^2  \frac {d\mu_{c}}{dc}|\leq K \|\eta\|_{L^2} \|v\|_{L^2}^2.$$

\medskip

\emph{Proof of \eqref{VIR2as1}.}
By the equation of $v$, we have
\begin{equation*}
\begin{split}
&-\frac 12 \frac d{dt} \int x v^2  = -\int x v_t v  
= -\int v_x \mathcal{L}_c (v x) +R_2 \quad \text{where}\\
&R_2= \int (f'(Q_c+\eta)-f'(Q_c)) x v_x v 
+\frac 12 (\rho'-c)\int  v^2 -c' \int x v Q_c -c' \int xv \eta.
\end{split}
\end{equation*}
First, by straightforward calculations, and using \eqref{decayQas1}, \eqref{surfas1}
$$-2 \int v_x \mathcal{L}_c (v x)  =    \int (3v_x^2+ v^2
-  f'(Q_c) v^2 -   xQ_c' f''(Q_c) v^2)
\geq \int (3 v_x^2+ v^2)- K \int v^2 e^{-\frac {\sqrt{c_0}}2|x|} .$$
Now, we estimate $R_2$:
\begin{equation*}\begin{split}
\left|\int (f'(Q_c+\eta)-f'(Q_c)) x v_x v \right|
&\leq K \|v\|_{L^\infty} \int |x\eta| |v_x|\leq K \|v\|_{H^1}^2 \|x\eta\|_{L^2}\\
&\leq  K \|v\|_{H^1}^2 \|x^2\eta\|_{L^2}^{\frac 12}\|\eta\|_{L^2}^{\frac 12}
\leq \frac 1{10} \|v\|_{H^1}^2,
\end{split}\end{equation*}
for $\alpha_0$ small enough, using \eqref{etapetit} and
 \eqref{decayvNL} (the constant in \eqref{decayvNL} does
not depend on $\alpha_0$).
We also have for $\alpha_0$ small, from \eqref{etapetit} and \eqref{orthovNL},
$ 
\frac 12 |\rho'-c|\int  v^2
\leq \frac 1{10} \int v^2.
$
From \eqref{orthovNL} and \eqref{decayQas1},
$
\left|c' \int x v Q_c\right|\leq 
\frac 1{10} \int v^2 + K \int v^2 e^{-\frac {\sqrt{c_0}}2 |x|}. 
$
Next,
$
|c'\int x v \eta|\leq K \|v\|_{L^2} \int |x\eta | |v|
\leq \frac 1{10} \|v\|_{H^1}^2$ is controled as above.
In conclusion, we have proved:
$$-\frac 12 \frac d{dt} \int x v^2 \geq    \frac 12 \int (v_x^2+v^2)   - K_0 \int  v^2e^{-\frac {\sqrt{c_0}}2 |x|}.$$

Now, fix $B>0$ such that $K_0 e^{-\frac {\sqrt{c_0}}2 B} \leq \frac 14$. Then, we obtain 
$$-\frac 12 \frac d{dt} \int x v^2 \geq   \frac 14 \int (v_x^2+v^2)   - K_0 \int_{|x|<B} v^2.$$

\medskip

\noindent\emph{Proof of Theorem 2.}
Consider $u(t)$ as in Theorem 2 with $\alpha_0>0$, small enough so that,
for all $t\geq 0$, $c(t)\in [c_0-\sigma_0,c_0+\sigma_0]\subset (0,c_*)$,
for $\sigma_0>0$ small enough so that Claim \ref{CL1as1} and Lemma \ref{CL100asz1} apply to $c(t)$.
Let 
$$V(t)=-\frac 12 \int (\mu_c+\varepsilon_0 x) v^2,$$
$$\text{where}\quad
\varepsilon_0=\frac 1{2} {\lambda_3}^2 \inf \{\mu_{c}'(x); |x|<B, c\in [c_0-\sigma_0,c_0+\sigma_0]\} >0
.$$
Then, from Lemma \ref{CL100asz1} and the definition of $\varepsilon_0$ we have for all $t$,
\begin{equation*}\begin{split}
V'(t)&\geq \lambda_3 \int v^2 \mu_c' +\lambda_3 \varepsilon_0 \int (v_x^2+v^2)
-\frac 1{\lambda_3} \|v\|_{H^1}^2\|\eta\|_{L^2}-\frac 1{\lambda_3} \varepsilon_0 \int_{|x|\leq B} v^2
\\ & \geq \lambda_3 \varepsilon_0 \int (v_x^2+v^2)
-\frac 1{\lambda_3} \|v\|_{H^1}^2\|\eta\|_{L^2}.
\end{split}\end{equation*}
Now, we choose $\alpha_0>0$ small enough so that by \eqref{etapetit},
$\frac 1{\lambda_3}\|\eta(t)\|_{L^2}\leq \frac 12 \lambda_3 \varepsilon_0$. Thus,
\begin{equation}\label{jfinal}
V'(t)\geq  \varepsilon_1 \int (v_x^2+v^2),
\quad \varepsilon_1=\frac 12 \lambda_5\varepsilon_0.
\end{equation}
By \eqref{decayvNL}, $V(t)$ is uniformly bounded on $\mathbb{R}$,
$\lim_{t\to \pm \infty} V(t)=V_\pm \infty$ and 
thus
\begin{equation}\label{uti}
\int_{-\infty}^{+\infty} \int (v_x^2+v^2)\leq \frac 1{\varepsilon_1}
(V_+ - V_-).
\end{equation}
Thus, there exist $t_n\to +\infty$ such that 
$v(t_n)\to 0,$ as $n\to +\infty$
in $H^1(\mathbb{R})$ and from  \eqref{decayvNL},
$V_+=\lim_{n\to + \infty} V(t_n)=0$. Similarly, $V_-=0$. 
Using \eqref{uti} again, we obtain
$$\forall t,x\in \mathbb{R},\quad v(t,x)\equiv 0.$$
From \eqref{orthovNL},
$ \forall t\in \mathbb{R},$ $
\eta(t)= 0,$ $ c'(t)=0,$ $ \rho'(t)=c(t).$
Thus, by \eqref{decomp}, $u(t,x)=Q_{c(0)}(x-c(0)t-\rho(0))$ is a soliton solution. This concludes the proof of Theorem 2.

\section{Asymptotic stability - Proof of Theorem 1}
The proof of the asymptotic stability is based on the nonlinear Liouville property
as in \cite{MMarchives}.

For a general nonlinearity, we do not use the direct approach used
 in  \cite{MMnonlinearity}. Indeed, for this approach, we would need spectral 
information on an linear operator related to $\mathcal{L}$, which we are not able to prove  in general. In contrast, the dual problem introduced in Section 3
can be understood for general nonlinearity, since the underlying linear operator is always nonnegative
(see Lemma \ref{VIRLIN}). This is an intrisic property of the dual problem.

Since working with the dual problem requires more regularity on the solution, 
we cannot work directly on the original $H^1$ solution.
Thus the proof of Theorem 1 consists in using Theorem 2 on limiting objects,
which are more regular than the solution itself.

However, we point out that the proof presented here is simpler than the one
in \cite{MMarchives}. Indeed, the convergence of $u(t_n)$ to an asymptotic object $\tilde u(t)$ is obtained
by monotonicity properties (such as Lemma \ref{LE3as1}) 
and not by the arguments of well-posedness
for the Cauchy problem for \eqref{fkdv} in $H^s$ ($0<s<1$)
and localization as in \cite{MMarchives}.

\medbreak 

We claim the following

\begin{proposition}[Convergence to a compact solution]\label{PROP3as1}
Under the assumptions of Theorem 1,
for any sequence $t_n\to +\infty$,
there exists  a subsequence $(t_{\phi(n)})$ and
$\tilde u_0\in H^1(\mathbb{R})$ such that for all $A>0$,
\begin{equation}
u(t_{\phi(n)},x+\rho(t_{\phi(n)})) \to \tilde u_0 
\quad \text{in $H^1(x>-A)$ as $n\to +\infty$},
\end{equation}
where
$c(t)$, $\rho(t)$ are associated to the decomposition of $u(t)$ as
in Lemma \ref{DECOMPas1}.

Moreover, the solution $\tilde u(t)$ of \eqref{fkdv} corresponding to  $\tilde u(0)=\tilde u_0$ is global
$(t\in\mathbb{R})$ and there exists $K>0$ such that 
\begin{equation}\label{proputilde}\begin{split}
\forall t \in \mathbb{R},\quad &
\|\tilde u(t,.+\tilde \rho(t))-Q_{c_0}\|_{H^1}\leq \alpha_0,\\
\forall t, x\in \mathbb{R},\quad &|\tilde u(t,x+\tilde \rho(t))|\leq K \exp\left(-\tfrac {\sqrt{c_0}}{16}
 |x|\right)
 \end{split}\end{equation}
 where 
 $\tilde c(t)$, $\tilde \rho(t)$ are associated to the decomposition of $\tilde u(t)$ as in Lemma \ref{DECOMPas1} and  $\tilde \rho(0)=0$.
\end{proposition}

Let us first prove Theorem 1 assuming Proposition \ref{PROP3as1} and then 
prove Proposition \ref{PROP3as1}.

\medskip

\noindent\emph{Proof of Theorem 1 assuming Proposition \ref{PROP3as1}.}
Let  $u(t)$  satisfy  the assumptions of Theorem 1 and 
 $\alpha_0>0$ small enough so that Theorem 2 holds.

\smallskip

From Proposition \ref{PROP3as1}, for any sequence $t_n\to +\infty$ 
there exists a subsequence $t_{n'}$ and $\tilde c_0$ such that $c(t_{n'})\to \tilde c_0$, and $\tilde u_0\in H^1(\mathbb{R})$ such that
$u(t_{n'},.+\rho(t_{n'}))-\tilde u_0\to 0$ in $H^1(x>-A)$, for any $A>0$.
Moreover, the solution $\tilde u(t)$ associated to $\tilde u(0)=\tilde u_0$
satisfies \eqref{proputilde} and $\tilde c(0)=\tilde c_0$,
$\tilde \rho(0)=0$.

Now we apply Theorem 2 to the solution $\tilde u(t)$. It follows that
$\tilde u(t)=Q_{c_1}(x-x_1-c_1 t)$. By uniqueness of the decomposition in Lemma \ref{DECOMPas1} applied to $\tilde u(0)$, we have
$c_1=\tilde c_0$ and $x_1=0$.

Therefore, 
$u(t_{n'},.+\rho(t_{n'}))-Q_{\tilde c_0}\to 0$ in $H^1(x>-A)$, for any $A>0$,
or equivalently, $u(t_{n'},.+\rho(t_{n'}))-Q_{c(t_{n'})}\to 0$ in $H^1(x>-A)$ for any $A>0$.
Thus, this being true for any sequence $t_n\to +\infty$, 
it follows that, for any $A>0$,
\begin{equation*}
u(t,.+\rho(t))-Q_{c(t)}\to 0\quad \text{in $H^1(x>-A)$ as $t\to +\infty$.}
\end{equation*}

\smallskip

Now, we observe that $\int Q_{c(t)}^2\to M_+>0$ as $t\to +\infty$.
This follows from   monotonicity arguments. See proof of Proposition 3
in \cite{MMarchives} and also step 3 of the proof of Proposition \ref{PROP3as1}.

Assuming now that there exists $\sigma_0>0$ such that
$c\mapsto \int Q_c^2$ is not constant in any interval
$I\subset [c_0-\sigma_0,c_0+\sigma_0]$. 
By possibly taking a smaller $\alpha_0>0$
so that $c(t)\in [c_0-\sigma_0,c_0+\sigma_0]$ for all $t$,
it follows from the continuity of $c(t)$ 
 that $c(t)$ has a limit as $t\to +\infty$.  

\smallskip

Finally, using the arguments of the proof of Proposition 3 in \cite{MMarchives},
we improve the convergence result to finish the proof of Theorem 1.

%\medskip

%

%\noindent\emph{Remark.} The conclusion of Proposition \ref{PROP3as1} holds under the only assumption that
%for $c_0>0$, a positive solution $Q_{c_0}$ of \eqref{ellipticas1} with $c=c_0$ exists.
%Indeed, \eqref{equivbis} is not used in the proof.

\medskip

\noindent\emph{Proof of Proposition \ref{PROP3as1}.} 
We consider a solution $u(t)$ satisfying the assumptions of Theorem 1.
First, we apply Lemma \ref{DECOMPas1} to $u(t)$: there exists $c(t)$, $\rho(t)$
satisfying \eqref{decomptt}--\eqref{decompttbis}, in particular,
there exists $K>0$ such that
\begin{equation}
\label{prin1}
\forall t\geq 0, \quad \|u(t,.+\rho(t))-Q_{c_0}\|_{H^1} \leq K \alpha_0.
\end{equation}

Let $t_n\to +\infty$. The sequence
$u(t_n,.+\rho(t_n))$ being bounded in $H^1$, there exists a subsequence of $(t_n)$
(still denoted by $(t_n)$)
and $\tilde u_0\in H^1(\mathbb{R})$ such that 
$$u(t_n,.+\rho(t_n))\rightharpoonup \tilde u_0 \quad \hbox{in $H^1$ weak.}$$
Let $\tilde u(t)$ be the solution of \eqref{fkdv} corresponding to $\tilde u(0)=\tilde u_0$
and defined on the maximal time interval $(-{\tilde T_-},{\tilde T_+})$.

\medskip

\emph{Step 1.} Exponential decay and strong convergence in $L^2$ on the right.

Consider the function $\psi$ defined on $\mathbb{R}$ by 
\begin{equation}\label{defpsias1}
\psi(x)=\frac 2 {\pi} \arctan\left(\exp\left( \frac x4\right)\right),
\quad\text{so that ${\rm lim}_{+\infty} \psi=1$, ${\rm lim}_{-\infty} \psi=0$.}
\end{equation}
Following Step 2 of the proof of Theorem 1 in \cite{MMnonlinearity} and the  monotonicity arguments
for \eqref{fkdv} (see Lemma \ref{LE3as1}), we have for all $x_0>0$,
\begin{equation}\label{limitas1}
\limsup_{t\to +\infty} \int(u_x^2+u^2)(t,x+\rho(t))\psi(\sqrt{{c_0}}(x-x_0))dx \leq K \, {\rm exp}\left(-\frac {\sqrt{{c_0}}}4{x_0}\right).
\end{equation}
Now, we prove the following

\begin{equation}\label{FORMexpas1}
\text{for all $A>0$,}\quad
u(t_n,.+\rho(t_n))\to \tilde u_0\quad \text{in $L^2(x>-A)$,}
\end{equation}
\begin{align}
 \forall t_0\in [0,\tilde T_+),\quad 
 &
\sup_{t\in [0,t_0]}\int (\tilde u_x^2+\tilde u^2)(t,x){\rm exp}\left(\frac {\sqrt{{c_0}}}4x\right) dx \leq K(t_0)<+\infty, \label{decayteras1} 
\\
&
\sup_{t\in [0,t_0]}\left\|\tilde u(t,x)  {\rm exp}\left(\tfrac {\sqrt{{c_0}}}8{x }\right)\right\|_{L^\infty} \leq K(t_0)<+\infty.\label{decaybisas1}
\end{align}

Proof of \eqref{FORMexpas1}.
Since $u(t_n,.+\rho(t_n))\rightharpoonup \tilde u_0$ in $H^1$ weak, 
we have
$u(t_n,.+\rho(t_n))\to \tilde u_0$ in $L^2_{loc}(\mathbb{R})$, and thus
by \eqref{limitas1}, we obtain, 
\begin{equation}\label{FORMas1}
\text{for all $A>0$,}\quad
u(t_n,.+\rho(t_n))\to \tilde u_0\quad \text{in $L^2(x>-A)$.}
\end{equation}

Proof of \eqref{decaybisas1}.  From \eqref{limitas1} and weak convergence in $H^1$, we have
for all $x_0>0$,
\begin{equation}\label{DDas1}
 \int (\tilde u_{0x}^2+\tilde u_0^2)(x)\psi(\sqrt{{c_0}}(x-x_0)) dx \leq
  K \, {\rm exp}\left(-\frac {\sqrt{{c_0}}}4{x_0}\right).
\end{equation}
Now, we prove a  similar estimate for $\tilde u(t)$, i.e. \eqref{decaybisas1} for 
$t\in [0,{\tilde T_+})$, with a rough constant and without using monotonicity arguments.
This kind of property is quite well-known for the gKdV equation (see Kato \cite{KATO}). 

Let $0<t_0<{\tilde T_+}$. Note that $\sup_{[0,t_0]}\| \tilde u(t)\|_{L^\infty}\leq 2 \sup_{[0,t_0]}\| \tilde u(t)\|_{H^1}\leq K(t_0)$,   and so 
$| \tilde uf(\tilde u)|+|F(\tilde u)|\leq K(t_0) \tilde u^2$. Thus, using   $0<\psi'<K \psi$
and $|\psi'''|\leq K \psi$,
by the computations of the proof of Lemma \ref{LE3as1}, we have, for all $x_0>0$,
\begin{align}
    &\frac 1{\sqrt{{c_0}}} \frac d {dt} \int \tilde u^2(t,x) \psi(\sqrt{{c_0}}(x-x_0))   
    =\int (-3 \tilde u_x^2 + 2 (\tilde uf(\tilde u)-F(\tilde u)))\psi'(\sqrt{{c_0}}(x-x_0))
    \nonumber
    \\ & \quad +{c_0} \int \tilde u^2 \psi'''(\sqrt{{c_0}}(x-x_0))  \leq K(t_0) \int \tilde u^2 \psi(\sqrt{{c_0}}(x-x_0)),\label{prin2}\\
    &\frac 1{\sqrt{{c_0}}} \frac d {dt} \int (\tilde u_x^2-2F(\tilde u)) \psi(\sqrt{{c_0}}(x-x_0)) 
       \nonumber\\ &\leq K(t_0) \int (\tilde u_x^2-2F(\tilde u))   \psi(\sqrt{{c_0}}(x-x_0))
      +K^2(t_0) \int \tilde u^2 \psi(\sqrt{{c_0}}(x-x_0)).\label{prin3}
\end{align}
First, we deduce from \eqref{DDas1} and \eqref{prin2} that 
$\forall t\in [0,t_0]$, $\forall x_0>0$,
$\int \tilde u^2(t,x) \psi(\sqrt{{c_0}}(x-x_0))  \leq K(t_0)\, {\rm exp}\left(-\frac {\sqrt{{c_0}}}4{x_0}\right) $. Then, by \eqref{prin3}, we obtain 
\begin{equation}\label{plustard}
\int (\tilde u_x^2+\tilde u^2)(t,x) \psi(\sqrt{{c_0}}(x-x_0))dx\leq K(t_0)\, {\rm exp}\left(-\frac {\sqrt{{c_0}}}4{x_0}\right).
\end{equation}
By \eqref{truc},  we have, for $\delta_1>0$, $\forall t\in [0,t_0]$, $\forall x_0>0$,
$$\int_{x<x_0} (\tilde u_x^2+\tilde u^2)(t,x){\rm exp}\left( \tfrac {\sqrt{{c_0}}}4{x}\right) dx \leq \tfrac 1{\delta_1} \exp\left( \tfrac {\sqrt{{c_0}}}4{x_0}\right)\int (\tilde u_x^2+\tilde u^2)(t,x) \psi(x-x_0)dx\leq \tfrac {1}{\delta_1} K(t_0),$$
and thus, passing to the limit $x_0\to +\infty$, \eqref{decayteras1} is proved.
Finally, by $\|w\|_{L^\infty(x>x_0)}^2\leq 2 \|w\|_{L^2(x>x_0)}\|w_x\|_{L^2(x>x_0)}$,
and \eqref{plustard} we also obtain
the   pointwise estimate \eqref{decaybisas1}.

\medskip

\emph{Step 2.} Strong convergence of $u(t_n+t,.+\rho(t_n))$ to $\tilde u(t)$ on the right.

\begin{lemma}\label{strongas1}
The solution $\tilde u(t)$ is global, i.e. ${\tilde T_-}={\tilde T_+}=+\infty$.
Moreover, for all $t\in \mathbb{R}$,
\begin{align*}
& \inf_{r\in \mathbb{R}} \|\tilde u(t,.+r)-Q_{c_0}\|_{H^1} \leq K \alpha_0,\\
&\text{ for all $A>0$,}\quad  u(t_n+t,.+\rho(t_n))\to \tilde u(t) \quad \text{in $H^1(x>-A)$
as $n\to +\infty$
,}\\
& \tilde \rho(0)=0,\quad  
\rho(t_n+t)-\rho(t_n)\to \tilde \rho(t)\quad \text{as $n\to +\infty$},
\end{align*}
where $\tilde \rho(t)$ is associated to the decomposition of $\tilde u(t)$ as in Lemma \ref{DECOMPas1}.
\end{lemma}

The proof of Lemma \ref{strongas1} contains the main new arguments.

\medskip

\noindent\emph{Proof of Lemma \ref{strongas1}.}
 For any $t\in (-\tilde T_-,\tilde T_+)$, we set
\begin{equation}\label{defvn}
v_n(t,x)=u(t_n+t,x+\rho(t_n))-\tilde u(t,x).
\end{equation}
Then, from the equation of $u(t)$ and $\tilde u(t)$ and \eqref{FORMexpas1},  $v_n(t)$ satisfies
\begin{align}
&  \partial_t v_{n}=-\partial_x (\partial_{x}^2 v_{n}+f(\tilde u+v_n)-f(\tilde u)),
\quad t\in (-\tilde T_-,\tilde T_+), x\in \mathbb{R}\label{prin5}\\
&  \int v_n^2(0) \psi(\sqrt{c_0} x) \to 0 \quad \text{as $n\to +\infty$.}\label{prin6}
\end{align}

\noindent\emph{Convergence in $L^2$ at the right for $t\geq 0$.}
Let  $0<t_0<{\tilde T_+}$.
We prove the following estimate:
\begin{equation}\label{prin7}
\sup_{t\in [0,t_0]} \int v_n^2(t) \psi(\sqrt{{c_0}}\, x)\leq K(t_0) \int v_n^2(0) \psi(\sqrt{{c_0}}\, x).
\end{equation}
Note that $\forall t\in [0,t_0]$, $\|\tilde u(t)\|_{L^\infty}
\leq K \|\tilde  u(t)\|_{H^1} \leq K(t_0)$, and since $f$ is $C^2$ and $f(0)=f'(0)=0$
$|f(\tilde u+v_n)-f(\tilde u)|\leq K |v_n|$ and $|F(v_n)|\leq K |v_n|^2$, and 
$|f(\tilde u+v_n)-f(\tilde u)-f(v_n)|\leq K |\tilde u| |v_n|$.

By   computations similar to the ones in  the
proof of Lemma \ref{LE3as1}, we have
\begin{equation*}\begin{split}
\frac 1{\sqrt{{c_0}}} \frac d{dt} \int v_n^2 \psi(\sqrt{{c_0}}\, x)   &=
-3 \int v_{nx}^2 \psi'(\sqrt{{c_0}}\, x)
+ {c_0} \int v_n^2 \psi'''(\sqrt{{c_0}}\, x) \\&+ 2\int (f(\tilde u{+}v_n)-f(\tilde u)) (v_n \psi(\sqrt{{c_0}}\, x))_x.
\end{split}\end{equation*}
We claim the following estimate of the nonlinear term:
\begin{equation}\label{prin8}
2\int (f(\tilde u{+}v_n)-f(\tilde u)) (v_n \psi(\sqrt{{c_0}}\, x))_x
\leq \frac 1{10}Ê\int v_{nx}^2  \psi'(\sqrt{{c_0}} \,x ) + K\int  v_n^2 \psi(\sqrt{{c_0}} \,x )
\end{equation}
Indeed, by direct computations;
\begin{align*}
& \int (f(\tilde u{+}v_n)-f(\tilde u)) (v_n \psi(\sqrt{{c_0}}\, x))_x      \\ &
=\int (f(\tilde u{+}v_n)-f(\tilde u)) (\sqrt{{c_0}} v_n \psi'(\sqrt{{c_0}}x)
+v_{nx} \psi(\sqrt{{c_0}}\, x)). \\
&= \sqrt{{c_0}} \int  ((f(\tilde u{+}v_n){-}f(\tilde u))v_n{-} F(v_n)) \psi'( \sqrt{{c_0}}  x)
+\int (f(\tilde u{+}v_n){-}f(\tilde u){-}f(v_n))  v_{nx} \psi(\sqrt{{c_0}}  x)  
 \\ & = \mathbf{I}+\mathbf{II}.  
\end{align*}
We have  
\begin{align*}
|\mathbf{I}|
& \leq
K(t_0)Ê\int v_n^2 \psi'(\sqrt{{c_0}}\, x) 
\\
  \left|
\mathbf{II}
\right| &\leq 
K(t_0) \int |\tilde u| |v_n| |v_{nx}| \psi(\sqrt{{c_0}} \,x )\\
& \leq K(t_0) \left\| \tilde u \sqrt{\frac { {\psi(\sqrt{{c_0}} \,x )}} {{ {\psi'(\sqrt{{c_0}}\, x)}}}} \right\|_{L^\infty}
\int |v_n| |v_{nx}| \sqrt{\psi(\sqrt{{c_0}} \,x )\psi'(\sqrt{{c_0}} \,x )}\\
&\leq  K(t_0) \left\| \tilde u \sqrt{\frac { {\psi(\sqrt{{c_0}} \,x )}} {{ {\psi'(\sqrt{{c_0}}\, x)}}}} \right\|_{L^\infty}
\left(\int v_{nx}^2  \psi'(\sqrt{{c_0}} \,x )\right)^{\frac 12} \left(\int  v_n^2 \psi(\sqrt{{c_0}} \,x )\right)^{\frac 12}.
\end{align*}
By the expression of $\psi$, we have
$\sqrt{\frac { {\psi(\sqrt{{c_0}} \,x )}} {{\psi'(\sqrt{{c_0}}\, x)}}}
\leq K (1+ e^{\frac {\sqrt{{c_0}}} 8 x})$
and thus from \eqref{decaybisas1}, we obtain 
\begin{align*}
  \left|
\mathbf{II}
\right|&\leq 
K(t_0) \left(\int v_{nx}^2  \psi'(\sqrt{{c_0}} \,x )\right)^{\frac 12} \left(\int  v_n^2 \psi(\sqrt{{c_0}} \,x )\right)^{\frac 12}
 \\& \leq \frac 1{10}Ê\int v_{nx}^2  \psi'(\sqrt{{c_0}} \,x ) + K(t_0)\int  v_n^2 \psi(\sqrt{{c_0}} \,x ) .
\end{align*}
Thus, \eqref{prin8} is proved and by $|\psi'''|\leq K \psi$, we obtain
\begin{align}\label{beepas1}
\frac 1{\sqrt{{c_0}}} \frac d{dt} \int v_n^2 \psi(\sqrt{{c_0}}\, x)    \leq 
-2 \int v_{nx}^2 \psi'(\sqrt{{c_0}}\, x)+K(t_0) \int v_n^2 \psi'(\sqrt{{c_0}}\, x).
\end{align}
This gives
\eqref{prin7}, 
and so, by \eqref{prin6}, we obtain
\begin{equation}\label{prin17}
 \sup_{t\in [0,t_0]} \int v_n^2(t) \psi(\sqrt{{c_0}}\, x)\to 0 \quad \text{as
$n\to +\infty$}.
\end{equation}
In particular,
for all $A>0$, as $n\to +\infty$,
$$u(t_n+t,.+\rho(t_n))\to \tilde u(t)  \text{ in $L^2(x>-A)$ and } 
u(t_n+t,.+\rho(t_n))\rightharpoonup \tilde u(t,.) \text{ in $H^1$ weak,}
$$
by the uniform $H^1$ bound on  $u(t)$, and by \eqref{prin1},
$$\forall t\in [0,\tilde T_+),\quad
\inf_{r\in \mathbb{R}} \|\tilde u(t,.+r)- Q_{c_0}\|_{H^1}\leq \alpha_0,\quad \text{
and so ${\tilde T_+}=+\infty$}.$$

\noindent\emph{Convergence in $L^2$ at the right for $t\leq 0$.}
Let $-{\tilde T_-}<t_1<0$. There exist $\tilde u_1(0)\in H^1$ and a subsequence $(t_{\phi(n)})$
such that
$$u(t_{\phi(n)}+t_1,.+\rho(t_{\phi(n)})\rightharpoonup \tilde  u_1(0) 
\quad \text{in $H^1$ weak as $n\to +\infty$}.$$
We reproduce on $\tilde u_1(0)$ the analysis done so far on $\tilde u(0)$. In particular,
let $\tilde u_1(t)$ be the solution of \eqref{fkdv} corresponding to $\tilde u_1(0)$
defined on $(-\tilde T_{1-},\tilde T_{1+})$.
It follows that $\tilde T_{1+}=+\infty$ and 
$$u(t_{\phi(n)},.+\rho(t_{\phi(n)}))\rightharpoonup \tilde  u_1(-t_1) \quad \text{in $H^1$ weak as $n\to +\infty$, and thus
$\tilde u_0=\tilde  u_1(-t_1).$}$$
By uniqueness of the $H^1$ solution of \eqref{fkdv}, we obtain
$\tilde  u_1(0)=\tilde u(t_1)$ and
 $u(t_{\phi(n)}+t_1,.+\rho(t_{\phi(n)}))\rightharpoonup \tilde u(t_1)$.
In fact, the convergence
$u(t_{n}+t_1,.+\rho(t_{n}))\rightharpoonup \tilde u(t_1)$ holds actually
for the whole sequence $(t_n)$.

As before, we obtain
$$\forall t>-\tilde T_-,\quad
\inf_{r\in \mathbb{R}} \|\tilde u(t,.+r)- Q_{c_0}\|_{H^1}\leq \alpha_0
\quad \text{and so $\tilde T_-=+\infty$}.$$
Therefore, we are able to  
  define $\tilde c(t)$, $\tilde \rho(t)$, associated to the decomposition of $\tilde u(t)$
as in Lemma \ref{DECOMPas1}. By continuity and uniqueness of the decomposition in $H^1$, we have
\begin{equation}\label{machh}
\tilde \rho(0)=0 \quad \text{and for all
$t\in \mathbb{R}$, } \rho(t_n+t)-\rho(t)\to \tilde\rho(t) \text{ as $n\to +\infty$}.
\end{equation}
In conclusion, in addition to \eqref{machh}, we have obtained so far, for all $t\in \mathbb{R}$,
\begin{align*}
& \inf_{r\in \mathbb{R}} \|\tilde u(t,.+r)-Q_{c_0}\|_{H^1} \leq K \alpha_0,\\
&\text{ for all $A>0$,}\quad  u(t_n+t,.+\rho(t_n))\to \tilde u(t) \quad \text{in $L^2(x>-A)$
as $n\to +\infty$
}.
\end{align*}

\noindent\emph{Convergence in $H^1$ at the right.} 
From the weak convergence  and \eqref{limitas1}, there exists $K>0$ such that 
\begin{equation}\label{prin10}
\forall x_0>0,\forall t\in \mathbb{R},\quad \int (\tilde u_x^2+\tilde u^2)(t,x+\tilde \rho(t)) \psi(\sqrt{c_0}(x-x_0)) dx 
\leq K \exp\left(-\tfrac {\sqrt{c_0}}{4} x_0\right), 
\end{equation}
and thus, as before,
\begin{equation}\label{prin14}
\forall t\in \mathbb{R},\forall x>0,\quad
|\tilde u(t,x+\tilde \rho(t))|\leq K \exp\left(-\tfrac {\sqrt{c_0}}{8} x\right).
\end{equation}

Let $v_n(t)$ be defined in \eqref{defvn} for all $t\in \mathbb{R}$.
We claim that  
\begin{equation}\label{prin16}
\forall  t\in \mathbb{R},\quad 
\int v_{nx}^2(t) \psi(\sqrt{{c_0}}x)dx\to 0\quad \text{as
$n\to +\infty$.}
\end{equation}
In particular, by \eqref{prin16} and \eqref{prin17}, this implies that 
for all $t\in \mathbb{R}$,
for all $A>0$,
$u(t_n+t,.+\rho(t_n))\to \tilde u(t)$ in $H^1(x>-A)$ and Lemma \ref{strongas1} follows.

Now, let us prove \eqref{prin16}. Let $t_0\in \mathbb{R}$.
It follows from  \eqref{beepas1} integrated on $[t_0-1,t_0]$ and  
\eqref{prin17}, that
$$\int_{t_0-1}^{t_0} \int v_{xn}^2(t,x)\psi'(\sqrt{{c_0}} x) dx dt 
\to 0\quad \text{as $n\to +\infty$.}$$
Thus, by \eqref{limitas1}, we obtain:
\begin{equation}\label{FORM3as1}
\int_{t_0-1}^{t_0} \int v_{xn}^2(t,x)\psi(\sqrt{{c_0}} x) dx dt 
\to 0\quad \text{as $n\to +\infty$.}
\end{equation}
Now, we claim for any $t_0-1\leq t\leq t_0$:
\begin{equation}\label{cla}\begin{split}
\int v_{nx}^2(t_0) \psi(\sqrt{{c_0}}x)dx
&\leq \int v_{nx}^2(t) \psi(\sqrt{{c_0}}x)dx
+ K(t_0) \int_{t_0-1}^{t_0} \int v_{nx}^2 (t') \psi(\sqrt{{c_0}}x)dxdt'\\ &
+ K(t_0) \sup_{t'\in [t_0-1,t_0]}\int v_n^2(t') \psi(\sqrt{{c_0}}x)dx.
\end{split}
\end{equation}
By \eqref{prin17} and \eqref{FORM3as1}, we find $\int v_{nx}^2(t_0) \psi(\sqrt{{c_0}}x)dx
\leq \int v_{nx}^2(t) \psi(\sqrt{{c_0}}x)dx+o(1)$, and thus
integrating on $t\in [t_0-1,t_0],$ using \eqref{FORM3as1} again, we prove 
 \eqref{prin16}.

\medskip

Now, let us prove \eqref{cla}.
Define
$$J(t)= \int (v_{nx}^2-2F(v_n))(t) \psi(\sqrt{{c_0}} x) dx,$$
so that 
\begin{align*}
\label{}
   &\frac 1 {\sqrt{{c_0}}}
   \frac  d{dt} J   = -3 \int v_{nxx}^2 \psi'(\sqrt{{c_0}} x) + {c_0}\int v_{nx}^2 \psi'''(\sqrt{{c_0}}x)
   \\ &+ 2 \int (f(\tilde u+v_n)-f(\tilde u))_x (v_{nx} \psi(\sqrt{{c_0}} x))_x
   - 2 \int (v_{nxx} + f(\tilde u +v_n)-f(\tilde u))(f(v_n)\psi(\sqrt{{c_0}} x))_x   \\
    &  \leq -2 \int v_{nxx}^2 \psi'(\sqrt{{c_0}} x) + K \int (v_{nx}^2 + v_n^2) \psi(\sqrt{{c_0}} x)
    \\ &+ 2 \int (f(\tilde u +v_n) -f(\tilde u)-f(v_n))_x v_{nxx} \psi(\sqrt{{c_0}} x),
\end{align*}
by controlling terms as in the proof of \eqref{prin8}
($\|v_n \sqrt{\psi(\sqrt{{c_0}} x)}\|_{L^\infty}^2
\leq K \int (v_{nx}^2+v_n^2) \psi(\sqrt{{c_0}} x)$). 
Now, we control the last term:
\begin{align*}
\label{}
    &	\int (f(\tilde u +v_n) -f(\tilde u)-f(v_n))_x v_{nxx} \psi(\sqrt{{c_0}} x)   \\
    &  = \int \tilde u_x (f'(\tilde u+v_n)-f'(\tilde u)) v_{nxx} \psi(\sqrt{{c_0}} x)
    + v_{nx} (f'(\tilde u+ v_n)-f'(v_n)) v_{nxx} \psi(\sqrt{{c_0}} x)\\
    & \leq K \int (|\tilde u_x| |v_n| + |\tilde u| |v_{nx}|) |v_{nxx}| \psi(\sqrt{{c_0}} x)\\
    & \leq \int v_{nxx}^2 \psi'(\sqrt{{c_0}} x) +
    K\left\|\tilde u \sqrt{\frac {\psi(\sqrt{{c_0}} x)}{\psi'(\sqrt{{c_0}} x)}}\right\|_{L^\infty}^2\int v_{nx}^2 \psi(\sqrt{{c_0}} x)\\ & +
    K\left\|v_n \sqrt{\psi(\sqrt{{c_0}} x)}\right\|_{L^\infty}^2\int ({\tilde u}_{x})^2 
    \frac {\psi(\sqrt{{c_0}} x)}{\psi'(\sqrt{{c_0}} x)}.
\end{align*}
We have $\|v_n \sqrt{\psi(\sqrt{{c_0}} x)}\|_{L^\infty}^2
\leq K \int (v_{nx}^2+v_n^2) \psi(\sqrt{{c_0}} x)$ and 
$ \frac { {\psi(\sqrt{{c_0}} \,x )}} {{\psi'(\sqrt{{c_0}}\, x)}}
\leq K (1+ e^{\frac {\sqrt{{c_0}}} 4 x})$. Thus,
 using \eqref{decayteras1}--\eqref{decaybisas1}, 
we obtain $$\frac d{dt} J(t)\leq K \int (v_{nx}^2+ v_n^2) \psi(\sqrt{c_0} x).$$
Integrating between $t$ and $t_0$ and using
$\int F(v_n)\psi(\sqrt{c_0} x)\leq K \int v_n^2 \psi(\sqrt{c_0} x)$,  
\eqref{cla} is proved. Thus, Lemma \ref{strongas1} is proved.

\medskip

\emph{Step 3.} Exponential decay of $\tilde u(t,x)$.
We prove
\begin{equation}\label{prin9}
\forall t,x\in \mathbb{R},\quad
|\tilde u(t,x+\tilde \rho(t))|\leq K \exp\left(-\tfrac {\sqrt{c_0}}{16} |x|\right).
\end{equation}
We claim
\begin{equation}\label{prin11}
\forall x_0>0,\forall t\in \mathbb{R},\quad \int \tilde u^2(t,x+\tilde \rho(t)) (1-\psi(\sqrt{c_0}(x+x_0))) dx 
\leq K \exp\left(-\tfrac {\sqrt{c_0}}{4} x_0\right).  
\end{equation}
Note that \eqref{prin9} is a direct consequence of \eqref{prin14}, \eqref{prin11} and the
global $H^1$ bound on $\tilde u(t)$ using $\|w\|_{L^\infty(x>x_0)}^2\leq 2 \|w\|_{L^2(x>x_0)}\|w_x\|_{L^2(x>x_0)}$.

\medskip

Proof of \eqref{prin11}.
Estimate \eqref{prin11}
was already proved in the same context in \cite{MMarchives} and \cite{MMT} (see for example \cite{MMT}, Lemma 7).
Let us sketch a proof.

We use  monotonicity arguments similar to the ones in Lemma \ref{LE3as1}.
Let $m_0=\int \tilde u_0^2$.
Let $x_0>0$ and $t_0\in \mathbb{R}$.
By $L^2$ norm conservation and  Lemma \ref{strongas1}, for $n(x_0)>0$ large enough, we have
\begin{align*}
\label{}
  & m_0  - \int \tilde u^2(t_0) (1-\psi(\sqrt{c_0} (x-\tilde \rho(t_0)+x_0))) 
  =  \int  \tilde u^2(t_0)  \psi(\sqrt{c_0} (x-\tilde \rho(t_0)+x_0))\\
  & \geq  \int u^2(t_n+t_0) \psi(\sqrt{c_0} (x-  \rho(t_0+t_n)+x_0))
   -   \exp\left(-\tfrac {\sqrt{c_0}} 4 x_0\right).
\end{align*}
By monotonicity properties on $u(t)$, for $n'\geq n$ so that $t_{n'}\geq t_n+t_0$, it follows that
\begin{align*}
\label{}
  & m_0  - \int \tilde u^2(t_0) (1-\psi(\sqrt{c_0} (x-\tilde \rho(t_0)+x_0))) 
    \\ &\geq  \int u^2(t_{n'} ) \psi(\sqrt{c_0} (x-\rho(t_{n'})+x_0+\tfrac {c_0}4 (t_{n'}-(t_n+t_0))))
   - K \exp\left(-\tfrac {\sqrt{c_0}} 4 x_0\right).
\end{align*}
Again from the convergence of $u(t_{n'},.+\rho(t_{n'})$ to $\tilde u(0)$, for
 $n'=n'(n,x_0)$ large enough, we have 
 $\int u^2(t_{n'}) \psi(\sqrt{c_0} (x-\rho(t_{n'})+x_0+\tfrac {c_0}4 (t_{n'}-(t_n+t_0))))
\geq m_0 -\exp\left(-\tfrac {\sqrt{c_0}} 4 x_0\right).$ This proves that
$
\int \tilde u^2(t_0) (1-\psi(\sqrt{c_0} (x-\tilde \rho(t_0)+x_0)))
\leq K \exp\left(-\tfrac {\sqrt{c_0}} 4 x_0\right),
$ thus \eqref{prin11} is proved.

\medskip

From Lemma \ref{strongas1} and \eqref{prin9}, Proposition \ref{PROP3as1} is proved.

\section{Multi-soliton case}

Now, we give a application of our results to the case of solutions containing
several solitons. Let  $N\geq 1$,
  $x_1,\ldots, x_N\in \mathbb{R}$, and 
\begin{equation}\label{stabxbx}
0<c_N^0<\ldots<c_1^0<c_*(f),\quad \forall j,~
\frac {\partial}{\partial c} \int {Q_c^2}_{|c=c_j} >0,
\end{equation}
it was proved in \cite{Martel} that there exists a unique solution $U(t)$
 in $H^1$ of  \eqref{fkdv} such that
\begin{equation}\label{defU}
\Big\|U(t)-\sum_{j=1}^N Q_{c_j}(.-c_j t - x_j) \Big\|_{H^1(\mathbb{R})} \to 0
\quad \text{as $t\to +\infty$.}
\end{equation}
This solution $U(t)$ is called a multi-soliton solution (in \cite{Martel}, the result
is proved only for the power case $f(u)=u^p$ for $p=2,3,4,5$
but the proof is exactly the same for a general $f(u)$ with stable solitons in the sense \eqref{stabxbx}).

The stability of such multi-soliton structures has been studied previously in
\cite{MMT}. Indeed, the main result in \cite{MMT} is that under assumption \eqref{stabxbx}, if 
\begin{equation}\label{th3bas22}
	\inf_{\substack{r_j\in \mathbb{R}\\ r_{j}-r_{j+1}>L_0}}
	\Big\|u(0)-\sum_{j=1}^N Q_{c_j^0}(.-r_j)\Big\|_{H^1} < \alpha_0,
	\end{equation}
	for $L_0$ large enough and $\alpha_0$ small enough,
then the solution $u(t)$ of \eqref{fkdv} satisfies
\begin{equation}\label{th3bas44}
	\forall t\geq 0,\quad
	\inf_{\substack{r_j\in \mathbb{R}\\ r_{j}-r_{j+1}>L_0}}
	\Big\|u(t)-\sum_{j=1}^N Q_{c_j^0}(.-r_j)\Big\|_{H^1} < A (\alpha_0
	+e^{-\gamma t}).
\end{equation}
Again the proof of this result in \cite{MMT} was for the power  case ($p=2,3,4$), but the same proof applies 
to a general $f(u)$ under assumption \eqref{stabxbx}.

In \cite{MMT}, the asymptotic stability of such multi-soliton was also proved,
but the proof was restricted to $p=2,3$ and $4$, since it was based on
\cite{MMnonlinearity} (linear Liouville argument).
As a direct consequence of Theorem 2 and the proof of
Theorem 1, we now extend the asymptotic stability result by the following.

\begin{theorem}[Asymptotic stability of multi-soliton solution]\label{TH3as1}
Assume that $f$ is $C^3$ and satisfies \eqref{surfas1}.
Let $N\geq 1$ and $0<c_N^0<\ldots<c_1^0<c_*(f)$.
There exist $L_0>0$ and $\alpha_0>0$ such that if
 $u(t)$ is a global  $(t\geq 0)$ $H^1$ solution of \eqref{fkdv}
satisfying \eqref{th3bas44}
then the following hold.
\begin{enumerate}
\item Asymptotic stability in the energy space.
There exist  $t\mapsto c_j(t)\in (0,c_*(f))$, $t\mapsto \rho_j(t)\in \mathbb{R}$ such that
\begin{equation}
\label{th3-1}
u(t)-\sum_{j=1}^NQ_{c_j(t)}(.-\rho(t))\to 0\quad \text{in $H^1(x>\tfrac {c_N^0}{10} t)$ as $t\to +\infty$.}
\end{equation}
\item Convergence of the scaling parameter.
Assume further that there exists $\sigma_0>0$ such that
$c\mapsto \int Q_c^2$ is not constant in any interval
$I\subset [c_j-\sigma_0,c_j+\sigma_0]$. 
Then, by possibly taking a smaller $\alpha_0>0$, 
there exits $c_{j,+}\in (0,c_*(f))$ such that $c(t)\to c_{j,+}$ as $t\to +\infty$.
\end{enumerate}
\end{theorem}

\noindent\emph{Sketch of the proof.}
The proof of Theorem \ref{TH3as1} does not use any new argument with respect to
Theorems \ref{TH1as1} and \ref{PROP2as1} and the proof of the main results
in \cite{MMT}.

The first observation is that assuming \eqref{th3bas44}, we have the analogue
of Lemma \ref{DECOMPas1}:
there exist $c_j(t)>0$, $\rho_j(t)\in C^1([0,+\infty))$ such that
\begin{equation}
\label{decomptt3}
\eta(t,x)=u(t,x)-\sum_{j=1}^N Q_{c_j(t)}(x-\rho_j(t)),
\end{equation}
satisfies, for all $t\in [0,T_0],$ for all $j=1,\ldots,N$,
\begin{align}
\label{propeps3}
    &   \int \tilde \chi_{c_j(t)}(x-\rho_j(t)) \eta(t,x)dx=
    \int Q'_{c(t)}(x-\rho_j(t)) \eta(t,x)dx=0, \\
    &	|c_j(t)-c_{j}^0|+\|\eta(t)\|_{H^1}\leq K_0 \alpha_0,\quad
    \rho_j(t)-\rho_{j+1}(t)>\frac {L_0}2+\sigma t
    ~(\sigma>0),
    \\
    &  |c'_j(t)|+|\rho'_j(t)-c_j(t)|\leq K_0 \left(\int \eta^2(t,x) e^{-|x-\rho_j(t)|}dx\right)^{\frac 12}. \label{decompttbis3}
\end{align}
Now, we prove asymptotic stability by considering various regions related to
the position of the solitons.

\smallskip

(a) Asymptotic stability around the first soliton on the right.

Here, we follow exactly the proof of Proposition \ref{PROP3as1}. 
Let $t_n\to +\infty$, for a subsequence $t_{\phi(n)}$, 
$ 
u(t_{\phi(n)},.+\rho_{1}(t_{\phi(n)}) \to \tilde u_{0,1},
$ and $\tilde u_1(t)$ solution of \eqref{fkdv} corresponding to 
$\tilde u_1(0)=\tilde u_{0,1}$ satisfies \eqref{proputilde}.
Indeed, in the proof of Proposition \ref{PROP3as1}, only the behavior of 
the solution $u(t)$ at the right of the soliton $Q_{c_1(t)}$ is concerned, the presence of
$N-1$ solitons on the left does not change the argument.
Thus, as in the proof of Theorem~\ref{TH1as1}, using Theorem \ref{PROP2as1},
we obtain $\tilde u_1(t)=Q_{c_{1,+}} (x-c_{1,+} t)$, where
$c_1(t_{\phi(n)})\to c_{1,+}$. Thus,
for any $A>0$,
$u(t,x+\rho_1(t))\to Q_{c_{1,+}}$ on $H^1(x>-A)$.
Finally, using only monotonicity arguments, we obtain
$$
u(t)- Q_{c_{1,+}}(.-\rho_1(t)) \to 0\quad \text{on } H^1(x>\tfrac 12 {(\rho_1(t)+\rho_2(t))})$$
see \cite{MMT}, Section 4.1 and \cite{MMnonlinearity}, proof of Theorem 1.

\smallskip

(b) Asymptotic stability on each solitons by iteration.
We prove the result on the other solitons by iteration on $j$ from $1$ to $N$
of the following statement:
\begin{equation}\label{induction}
\exists c_{j,+}\text { s.t.}\quad 
u(t)- Q_{c_{j,+}}(.-\rho_j(t)) \to 0\quad \text{on } H^1(x>\tfrac 12 {(\rho_j(t)+\rho_{j+1}(t))}),
\end{equation}
(if $j=N$, the convergence is on $H^1(x>\tfrac 1{10}  c_N^0 t))$.

Assume that \eqref{induction} holds for $1\leq j_0 <N$. Let us prove it for $j_0+1$.  
The only point that differs from the case of $j=1$ is the analogue
of Lemma \ref{strongas1} to prove strong convergence in $H^1$ on the right.

For any $t_n\to +\infty$,
there exists $c_{j_0+1}$, $\tilde u_{0,j_0+1}$ such that (up to a subsequence still denoted by
$t_n$):
$$
u(t_n,.+ \rho_{j_0+1}(t_n)) \to \tilde u_{0,j_0+1} \quad \text{in $L^2_{loc}$},
\quad c_j(t_n)\to c_{j_0+1},
$$
where $\tilde u_{0,j_0+1}$ has exponential decay on the right.
Set, for $j=1,\ldots,j_0$,
$$
R_{j}(t,x)=R_{j}^{n,j_0} (t,x)=
Q_{c_{j,+}}(x-c_{j,+} t - \rho_{j} (t_n) + \rho_{j_0+1}(t_n)),
$$
$$
v_n(t,x)=u(t_n+t,x+\rho_{j_0+1}(t_n))
-\tilde u_{j_0+1}(t,x) - \sum_{j=1}^{j_0} R_j(t,x),
$$
where $\tilde u_{j_0+1}=\tilde u$ is the solution  of \eqref{fkdv}
corresponding to $\tilde u_{0,j_0+1}$.

Following Proposition \ref{PROP3as1}, it is enough to prove
\begin{equation}\label{claimgg}
\int (v_{nx}^2 + v_n^2)(t,x) \psi(\sqrt{c_{j_0+1}} x )  dx\to 0 \quad 
\text{as $n\to +\infty$.}
\end{equation}
Proof of \eqref{claimgg}.
Let us just check convergence for $\int   v_n^2(t,x) \psi(\sqrt{c_{j_0+1}} x )  dx$, the case of $v_{nx}$ is treated as in Proposition \ref{PROP3as1}.
First, we have
$$
\partial_t v_n = -  \partial_x ( \partial^2 v_n + f(\tilde u + \sum_{j=1}^{j_0} R_j
+v_n) - f(\tilde u) - \sum_{j=1}^{j_0} f(R_j)), \quad \text{and}
$$
$$
\int v_n^2(0,x) \psi(\sqrt{c_{j_0+1}} x )  dx\to 0 \quad 
\text{as $n\to +\infty$.}
$$
Computing (energy method) $\frac d{dt} \int v_n^2(t,x) \psi(\sqrt{c_{j_0+1}} x )  dx$, as in the proof of Proposition \ref{PROP3as1}, the only term
which has to checked is:
\begin{equation*}\begin{split}
&\int (f(\tilde u + \sum_{j=1}^{j_0} R_j
+v_n) - f(\tilde u) - \sum_{j=1}^{j_0} f(R_j)) v_{nx} \psi(\sqrt{c_{j_0+1}} x )
=\\
& \int (f(\tilde u + \sum_{j=1}^{j_0} R_j
+v_n) - f(\tilde u) - f(\sum_{j=1}^{j_0} R_j+v_n)) v_{nx} \psi(\sqrt{c_{j_0+1}} x )
\\ &+
\int ( f(\sum_{j=1}^{j_0} R_j+v_n)
- \sum_{j=1}^{j_0} f(R_j) ) v_{nx} \psi(\sqrt{c_{j_0+1}} x )
=I + II.
\end{split}\end{equation*}
$$
|I|\leq C \int|\tilde u|(|v_n|+ \sum_{j=1}^{j_0}Ê|R_j|) |v_{nx}|\psi(\sqrt{c_{j_0+1}} x )
\leq 
C \int |\tilde u||v_n||v_{nx}|\psi(\sqrt{c_{j_0+1}} x ) + C e^{-\sigma (t_N+t)}.
$$
\begin{equation*}
\begin{split}
II & =
-\sqrt{c_{j_0+1}} \int (F(\sum_{j=1}^{j_0} R_j + v_n)  - F(\sum_{j=1}^{j_0} R_j) 
- v_n f(\sum_{j=1}^{j_0} R_j)) 
\psi'(\sqrt{c_{j_0+1}} x )\\
& + \int v_n (\sum_{j=1}^{j_0} f(R_j) - f(\sum_{j=1}^{j_0} R_j))\psi(\sqrt{c_{j_0+1}} x )\\ &
-\int \sum_{j=1}^{j_0} (R_{jx} (f(\sum_{j=1}^{j_0} R_j + v_n)- f(\sum_{j=1}^{j_0} R_j)) - v_n f'(\sum_{j=1}^{j_0}R_j)) \psi(\sqrt{c_{j_0+1}} x )
\\ & =II_1+II_2+II_3.
\end{split}
\end{equation*}
Then $
|II_1|\leq 
C \int v_n^2 \psi'(\sqrt{c_{j_0+1}} x ),
$ $ 
|II_2|\leq 
C e^{-\sigma (t_n+t)},
$ $ 
|II_3|\leq C\int v_n^2 \psi(\sqrt{c_{j_0+1}} x )
$  implies
$$
\frac d{dt} \int v_n^2(t,x) \psi(\sqrt{c_{j_0+1}} x )  dx
\leq C \int v_n^2(t,x) \psi(\sqrt{c_{j_0+1}} x )  dx + C e^{-\sigma (t_n+t)},
$$
and the conclusion.

\appendix

\section{Monotonicity results}
Define
$
\psi(x)=\frac 2 {\pi} \arctan(\exp(  x/4)),
$
so that ${\rm lim}_{+\infty} \psi=1$, ${\rm lim}_{-\infty} \psi=0$
and for all $x\in \mathbb{R}$, $\psi(-x)=1-\psi(x)$. 
Note also that by direct calculations
\begin{equation}\label{surphi} 
\psi'(x)=\frac 1{4 \pi {\rm cosh}(x/4)} >0,
\quad \psi'''(x)\le \frac {1} {16} \psi'(x),\end{equation}
\begin{equation}\label{truc}
\exists \delta_1>0,~\forall x<0,\quad 
\psi(x)\geq \delta_1 \exp\left(\frac x 4\right), \quad \psi'( x)\geq \delta_1 \exp\left(\frac x 4\right).
\end{equation}

\subsection{Monotonicity arguments on $u(t)$}

Let $u(t)$ be a solution of \eqref{fkdv} satisfying the assumptions of Lemma \ref{DECOMPas1}
for   $t\in [0,T_0]$.
Let   $x_0>0$. We define,  for $0\leq t\leq t_0\leq T_0$:
 $\psi_0(t,x)=\psi(\sqrt{c_0}(x-\rho(t_0)+\tfrac {c_0} 2(t_0-t)-x_0))$ and 
\begin{equation*}
I_{x_0,t_0}(t)
=\int u^2(t,x)\,\psi_0(t,x) dx,
\quad J_{x_0,t_0}(t)
=\int \bigl( u_x^2-2 F(u)+c_0 \,u^2\bigr)(t,x)\,\psi_0(t,x) dx.
\end{equation*}

\begin{lemma}\label{LE3as1}
There exists  $K=K(c_0)>0$ such that for $\alpha_0$ small enough,
 for all $0\leq t\leq t_0\leq T_0$,
\begin{equation}
I_{x_0,t_0}(t_0)-I_{x_0,t_0}(t) \le K \, {\rm exp}\left(-\frac {\sqrt{c_0}}4{x_0}\right), 
\quad J_{x_0,t_0}(t_0)-J_{x_0,t_0}(t) \le K \, {\rm exp}\left(-\frac {\sqrt{c_0}}4{x_0}\right). \label{le3deux}\end{equation}
\end{lemma}

\noindent\emph{Proof of Lemma \ref{LE3as1}.}
The proof is the same as the one of Lemma 3 in \cite{MMnonlinearity} for $u^p$,
we repeat it for a general nonlinearity $f(u)$.
By simple calculations, for $\phi:\mathbb{R}\rightarrow \mathbb{R}$
of class $C^3$, we have
\begin{align*}
\frac {d}{dt} \int u^2 \phi &= 
2\int u_t u \phi=-2\int (u_{xx}+f(u))_x u \phi
= 2 \int (u_{xx}+f(u)) (u_x \phi+u\phi') 
\\&= \int \Bigl(-3 u_x^2 + 2(uf(u)-F(u))\Bigr)\phi'+\int u^2\phi''',
\end{align*}
\begin{align*}
\frac d{dt} \int \Bigl(u_x^2-2F(u)\Bigr)\phi &= 
2\int \bigl(u_{xt}u_x-f(u)u_t\bigr)\phi=-2\int u_t\bigl(u_{xx}+f(u)\bigr)\phi
-2\int u_t u_x \phi' 
\\&=  -\int \bigl(u_{xx}+f(u)\bigr)^2 \phi' +2\int\bigl(u_{xx}+f(u)\bigr)_x u_x \phi'
\\&=  \int \Bigl[-\bigl(u_{xx}+f(u)\bigr)^2 -2 u_{xx}^2  
+2 u_x^2 f'(u)\Bigr] \phi'
+\int u_x^2 \phi'''.
\end{align*}
We obtain from the previous calculations and  (\ref{surphi}), for all $t\le t_0$,
\begin{align*}
\frac {d}{dt} \int u^2 \psi_0
& =-\int \left(3 u_x^2+\tfrac {c_0}2 u^2 -2(uf(u)-F(u)) \right)\psi_{0x}+\int u^2\psi_{0xxx}
\\ & \leq 
- \int  \left(3u_x^2+\frac {c_0}4 u^2
- 2(uf(u)-F(u))\right)\psi_{0x}.
\end{align*}
Let $R_0>0$ to be chosen later.

(i) For $t,x$ such that $|x-\rho(t)|\ge R_0$, 
by \eqref{decayQas1},  (\ref{decomptt}), 
\begin{equation*}
|u(t,x)|\le Q_{c_0}(x)+\|u(t)- Q_{c_0}\|_{L^\infty}\le Q_{c_0}(x)+\|u(t)- Q_{c_0}\|_{H^1}
\le  K( e^{-\frac  {\sqrt{c_0}}2 R_0}+ \alpha_0).
\end{equation*}
Therefore, for $\alpha_0$ small enough and $R_0$ large enough, we have,
for such $t,x$:
$ |uf(u)-F(u)|\le \frac {c_0}8 \, u^2$.
Now, $\alpha_0$ and $R_0$ are fixed to such values.

(ii) For $t,x$ such that $|x-\rho(t)|\le R_0$ then $|x-\rho(t_0)+\tfrac 12 (t_0-t) -x_0|\ge -|x-\rho(t)|+
|\rho(t)-\rho(t_0)+\sigma (t_0-t) -x_0|\ge -R_0+\frac {t_0-t}2+x_0$, and so
\begin{equation*}
|\psi_{0x}(t,x)|\le K e^{-\frac {\sqrt{c_0}}8(t_0-t)}e^{-\frac {\sqrt{c_0}}4 {x_0}}.
\end{equation*}
Therefore, since $\|u\|_{L^\infty}\le K$, we obtain
\begin{equation}\label{gyyy}
\frac {d}{dt} \int u^2 \psi_0
\leq -\int  \left(3u_x^2+\frac {c_0}8 u^2\right)\psi_{0x}
-K e^{-\frac {\sqrt{c_0}}8(t_0-t)}e^{-\frac {\sqrt{c_0}}4 {x_0}}
\leq -K e^{-\frac {\sqrt{c_0}}8(t_0-t)}e^{-\frac {\sqrt{c_0}}4 {x_0}}  .
\end{equation}
By integration between $t$ and $t_0$, we obtain \eqref{le3deux} for $\mathcal{I}_{x_0,t_0}$.

Similarly, using (\ref{surphi}), we have
\begin{align*}
\frac {d}{dt} \int \Bigl(u_x^2
-2F(u)\Bigr)\psi_0
&=  \int \Bigl[-\bigl(u_{xx}+f(u)\bigr)^2 -2 u_{xx}^2  
+2  u_x^2 f'(u)\Bigr]\psi_{0x}\\
& -\frac {c_0}2 \int \Bigl(u_x^2
-2 F(u)\Bigr)\psi_{0x}+\int u_x^2 \psi_{0xxx}\\
& \le - \int \Bigl[\bigl(u_{xx}+f(u)\bigr)^2 +2 u_{xx}^2  
+\frac {c_0}4 u_x^2-2 u_x^2 f'(u)
-{c_0}F(u)\Bigr] \psi_{0x}.
\end{align*}
Splitting in two regions ($|x-\rho(t)|\geq R_0$, $|x-\rho(t)|\leq R_0$) as before,
by the same argument, we control the nonlinear terms   so that by (\ref{gyyy})
\begin{equation*}
\frac {d}{dt} \int \Bigl( u_x^2-2F(u)+c_0\, u^2\Bigr)\psi_0
\le  -\int \left( 2 u_{xx}^2+\frac {c_0}8 u_x^2
+\frac {c_0^2}{16} u^2\right)\psi_{0x}
-K e^{-\frac {\sqrt{c_0}}8(t_0-t)}e^{-\frac {\sqrt{c_0}}4 {x_0}} .
\end{equation*}
Therefore, by integration, we obtain
\eqref{le3deux}.
Note for future use that for $0\leq t<t_0\leq T_0$,
\begin{align}
&\int  u^2(t_0) \psi_0(t_0)   + 
\frac 1{16} \int_{t}^{t_0} \int (u_x^2 + c_0 u^2)  \psi_{0x}   dt'\leq \int  u^2(t) \psi_0(t)   + K \exp\left(-\frac {\sqrt{c_0}}4 {x_0}\right)\label{futur1}\\
&\int (u_x^2+c_0 u^2)(t_0) \psi_0(t_0)   + 
\frac 1{16} \int_{t}^{t_0} \int (u_{xx}^2 + c_0 u_x^2 + c_0^2 u^2)  \psi_{0x}    dt'\label{futur2}	\\
&\leq \int (u_x^2(t)+u^2(t)) \psi_0(t)   + K \exp\left(-\frac {\sqrt{c_0}}4 {x_0}\right).
\label{futur3}
\end{align}

\subsection{Monotonicity arguments on the linearized problem. Proof of   \eqref{monov}}

Let $\eta(t)$ be as in Proposition \ref{PROPas1}. We  claim the following preliminary result.

\begin{claim}\label{esteta}
There exists $K>0$ such that, for all $t_0\in \mathbb{R}$,
$$
\int_{-\infty}^{t_0}\int (\eta_{xx}^2+c_0 \eta_x^2+c_0^2 \eta^2)
\exp\left(\frac {\sqrt{c_0}} 4 (x-\tfrac {c_0}2(t_0-t)\right) dxdt \leq   K.
$$
\end{claim}

\noindent\emph{Remark.} We obtain a gain of regularity on $\eta(t,x)$ using the decay assumption \eqref{dsds} and monotonicity arguments.

\medskip

\noindent\emph{Proof of Claim \ref{esteta}.}
Let $t_0\in \mathbb{R}$, $x_0>0$, and $\tilde x= \sqrt{{c_0}} (x-\frac {c_0}2 (t_0-t) -x_0)$.
Then, by similar calculations as in Lemma \ref{LE3as1}, using in particular
\eqref{surphi}, we have
\begin{align*}
     \frac d{dt} \int \eta^2 \psi(\tilde x) & \leq 
    - 3 \sqrt{{c_0}} \int  \eta_x^2 \psi'(\tilde x) -\frac 14{c_0^{3/2}}   \int \eta^2 \psi'(\tilde x) 
    \\&+\int \eta^2 (\sqrt{{c_0}} f'(Q_{c_0}) \psi'(\tilde x)
    - f''(Q_{c_0}) Q_{c_0}' \psi(\tilde x)) 
\end{align*}
\begin{equation}
\label{enerV}
    \frac d{dt} \int \eta_x^2 \psi(\tilde x) 
     \leq - 3 \sqrt{c_0} \int \eta_{xx}^2 \psi'(\tilde x) - \frac 14 c_0^{3/2} \int \eta_x^2 \psi'(\tilde x) - 2\int (f'(Q_{c_0}) \eta)_x (\eta_x \psi(\tilde x))_x.  
\end{equation}
By \eqref{surfas1} and \eqref{decayQas1}, we have
\begin{equation}\label{dekayp}
|f'(Q_{c_0}(x))|+|f''(Q_{c_0}(x)) Q_{c_0}'(x)|\leq K \exp(-\sqrt{c_0} (p-1)|x|).
\end{equation}
Using \eqref{dsds} and 
arguing as in \cite{yvanSIAM}, proof of Lemma 5, we obtain for $x_0>0$ 
(considering three regions),
\begin{equation*}
\frac d{dt} \int \eta^2 \psi(\tilde x)  \leq  -3  \sqrt{{c_0}} \int  \eta_x^2 \psi'(\tilde x)
-\frac 18{c_0^{3/2}}   \int \eta^2 \psi'(\tilde x) 
+K \exp\left(-\sqrt{{c_0}}\left(\frac {c_0}8 (t_0-t) -\frac {x_0}4\right)\right).
\end{equation*}
Integrating between $t<t_0$ and $t_0$, we obtain, for all $t$
\begin{equation*}\begin{split}
&\int \eta^2(t_0) \psi(\sqrt{{c_0}}(x-x_0)) +\frac 18\sqrt{c_0} \int_{t}^{t_0}
\int  (\eta_x^2+c_0\eta^2) \psi'(\tilde x)dt
\\&\leq \int \eta^2(t) \psi(\sqrt{{c_0}} (x-\tfrac {c_0}2 (t_0-t) -x_0))
+K \exp\left(- \frac {\sqrt{{c_0}}} 4{x_0}\right).
\end{split}\end{equation*}
Passing to the limit $t\to -\infty$, using \eqref{dsds} and then  \eqref{truc}, 
we find, for all $t_0$,
\begin{equation}\label{monoeta}
 \int_{-\infty}^{t_0} \int  (\eta_x^2+c_0 \eta^2) \psi'(\tilde x)dt
\leq K \exp\left( -\frac {\sqrt{{c_0}}} 4{x_0}\right).
\end{equation}
\begin{equation*}\begin{split}
&\int_{-\infty}^{t_0}\int_{x<x_0+\frac {c_0}2 (t_0-t)} \int (\eta_x^2+c_0 \eta^2) \exp\left(\frac {\sqrt{c_0}} 4 (x-\tfrac {c_0}2(t_0-t)\right) dxdt
\\ &\leq \frac 1{\delta_1} \exp\left( \frac {\sqrt{{c_0}}}4{x_0}\right) 
\int_{-\infty}^{t_0} \int(\eta_x^2+c_0 \eta^2) \psi'(\tilde x) dx dt\leq K.
\end{split}\end{equation*}
Let  $x_0\to +\infty$, we find,
for all $t_0\in \mathbb{R},$
\begin{equation}\label{expoetaint}
\int_{-\infty}^{t_0}\int (\eta_x^2+c_0 \eta^2) \exp\left(\frac {\sqrt{c_0}} 4 (x-\tfrac {c_0}2(t_0-t)\right) dxdt
\leq   K.
\end{equation}

Now, we use \eqref{enerV}.
We expand the nonlinear term as follows:
\begin{equation*}
\label{ }
\begin{split}
&2\int (f'(Q_{c_0}) \eta)_x (\eta_x \psi(\tilde x))_x =
2\int (f'(Q_{c_0})\eta_x + f''(Q_{c_0}) Q'_{c_0}\eta)(\eta_{xx} \psi(\tilde x) +\sqrt{c_0}\eta_x \psi'(\tilde x))\\
&=\int  \eta_x^2 (-f''(Q_{c_0}) Q_{c_0}'\psi(\tilde x) +\sqrt{c_0} f'(Q_{c_0})\psi'(\tilde x))
+2 \int \sqrt{c_0} f''(Q_{c_0}) Q'_{c_0} \psi'(\tilde x) \eta \eta_{x}
\\ & +2 \int f''(Q_{c_0}) Q'_{c_0} \psi(\tilde x) \eta \eta_{xx}
=\mathbf{I}+\mathbf{II}+\mathbf{III}.
\end{split}
\end{equation*}
Note that by \eqref{dekayp}, we have
\begin{equation}
\label{decayq}
(|f'(Q_{c_0}(x))|+|f''(Q_{c_0}(x)) Q_{c_0}'(x)|)\psi(\tilde x)\leq K \exp(-\sqrt{c_0} (p-1)|x|)
\psi(\tilde x) \leq K \psi'(\tilde x).
\end{equation}
Indeed, for $\tilde x\leq 0$, we have $\psi(\tilde x)\leq K \psi'(\tilde x)$ and for 
$\tilde x>0$, we have $0<\tilde x\leq x$ and so
 $\exp(-\sqrt{c_0} (p-1)|x|)\leq K \exp(-\sqrt{c_0}\tilde x)\leq K \psi'(\tilde x)$.

Thus, $\mathbf{I}+\mathbf{II} \leq   K \int (\eta_x^2+ \eta^2) \psi'(\tilde x)$
and
$$
\mathbf{III}\leq
\sqrt{c_0} \int \eta_{xx}^2 \psi'(\tilde x) + K \int  \eta^2  \psi'(\tilde x).
$$
From \eqref{enerV}, we obtain
\begin{equation}
\label{ouf}
    \frac d{dt} \int \eta_x^2 \psi(\tilde x) 
     + 2 \sqrt{c_0} \int \eta_{xx}^2 \psi'(\tilde x) 
     \leq K \int (\eta_x^2+ \eta^2) \psi'(\tilde x).  
\end{equation}
From \eqref{expoetaint}, there exists a sequence $t_n\to -\infty$ so that
$\int \eta_x^2(t_n) \psi(\tilde x) \to 0$ as $n\to -\infty$. Thus, integrating \eqref{ouf}
between $t_0$ and $t_n$ and passing to the limit as $n\to +\infty$, using \eqref{monoeta}, we
obtain, arguing as before, for all $t_0\in \mathbb{R}$, 
\begin{equation}\label{monoetatrois}
\int_{-\infty}^{t_0} \int \eta_{xx}^2(t) \psi'(\tilde x) dx dt 
\leq K \exp\left( -\frac {\sqrt{{c_0}}} 4{x_0}\right),
\end{equation}
\begin{equation}\label{monoetaquatre}
 \int_{-\infty}^{t_0}\int \eta_{xx}^2(t) \exp\left(\frac {\sqrt{c_0}} 4 (x-\tfrac {c_0}2(t_0-t))\right) dxdt
\leq   K.
\end{equation}
and Claim \ref{esteta} is proved.

\medskip

Assuming that $f$ is $C^3$, one can apply monotonicity arguments again on $\eta_{xx}(t)$
and conclude the result for $v$. If $f$ is only $C^2$, we use the equation of $v$.
Recall that we argue by density again.

\medskip

\noindent\emph{Proof of   \eqref{monov}.}
Setting $\tilde v(t)=v(t)-\alpha(t) Q_{c_0}$, we see that $\tilde v(t)$ satisfies
$$\tilde v_t=\mathcal{L}_{c_0} \tilde v_x.$$
First, by the definition of $\tilde v(t)$ and Claim \ref{esteta} (see \eqref{expoetaint}, 
\eqref{monoetatrois}),  we have
\begin{equation}\label{controlvun}
 \int_{-\infty}^{t_0} \int  \tilde v^2(t) \psi'(\tilde x) dx dt 
\leq K \exp\left( -\frac {\sqrt{{c_0}}} 4{x_0}\right),
\end{equation}
\begin{equation}\label{controlvdeux}
 \int \tilde v^2(t_n) \psi(\tilde x) dx \to 0 \quad
\end{equation}
for a sequence $t_n\to -\infty$, where $\tilde x$ is defined in Claim \ref{esteta}.
By the equation of $\tilde v$, we have as in the proof of Claim \ref{esteta}
$$
\frac d {dt} \int \tilde v^2 \psi(\tilde x) 
\leq - \frac 12  \sqrt{c_0} \int (\tilde v_x^2 +c_0 \tilde v^2)\psi'(\tilde x)
-\int \tilde v^2 (\sqrt{c_0} f'(Q_{c_0}) \psi'(\tilde x) + f''(Q_{c_0})Q'_{c_0}\psi(\tilde x))
$$
$$
\frac d {dt} \int \tilde v_x^2 \psi(\tilde x) 
\leq - \frac 12  \sqrt{c_0} \int (\tilde v_{xx}^2+ c_0 \tilde v_x^2) \psi'(\tilde x)
-\int \tilde v_x^2 (\sqrt{c_0} f'(Q_{c_0}) \psi'(\tilde x) - f''(Q_{c_0})Q'_{c_0}\psi(\tilde x)).
$$
Integrating on $(-\infty, t_0]$ and combining these estimates with \eqref{controlvun}, 
\eqref{controlvdeux}, arguing as in the proof of Claim \ref{esteta}, we obtain
for all $t\in \mathbb{R}$,
$$
\int (v_x^2+c_0 v^2)(t) \exp\left(\frac {\sqrt{c_0}} 4 x\right) dx\leq K.
$$
Using the transformation $x\to -x$, $t\to -t$
(the equation of $v$ and the assumptions are invariant by this transformation),
 \eqref{monov} is proved.

\subsection{Proof of (\ref{monovNL}).}
We are in the context of  Theorem 2. In particular,
we assume \eqref{stableas1} and \eqref{decayas1} on the solution $u(t)$.
Using the same arguments as in the proof of Claim \ref{esteta}, we first claim the following preliminary result of $u(t)$.

\begin{claim}\label{estutilde} There exists $K>0$ such that
\begin{equation}\label{agarder4}
\forall t\in \mathbb{R},\quad 
\int (u_{xx}^2+ c_0 u_x^2+c_0^2 u^2)(t)\exp\left(\frac {\sqrt{c_0}} 4 |x-\rho(t)|\right) dx
\leq K.
\end{equation}
\end{claim}

\noindent\emph{Proof of Claim \ref{estutilde}.} 
From \eqref{decayas1}, letting $t\to -\infty$ in \eqref{futur1},  we have
$$\int_{-\infty}^{t_0} \int (u_x^2 + c_0 u^2)  \psi_{0x}    dt
\leq  K \exp\left(-\frac {\sqrt{c_0}}4 {x_0}\right).$$
By \eqref{truc}, and then letting $x_0\to +\infty$,
%$$\int_{x<\rho(t_0)+x_0} u^2(t_0) \exp\left(\frac {\sqrt{c_0}} 4 (x-\rho(t_0))\right) dx 
%\leq \frac 1{\delta_1} \exp\left(\frac {\sqrt{c_0}} 4 x_0\right)\int u^2(t_0) \psi_{0}(t_0)
%\leq K$$ 
\begin{align*}
& \int_{-\infty}^{t_0} \int_{x<x_0+\rho(t_0)-\frac {c_0}2 (t_0-t)}
(u_x^2+c_0 u^2)(t) \exp\left(\frac {\sqrt{c_0}} 4 (x-\rho(t_0)+\tfrac {c_0} 2 (t_0-t))\right) dxdt	\\
& \leq \frac 1{\delta_1} \exp\left(\frac {\sqrt{c_0}} 4 x_0\right)
\int_{-\infty}^{t} \int (u_x^2+c_0 u^2) \psi_{0x} dxdt\leq K,
\end{align*}
\begin{equation}\label{agarder}
\int_{-\infty}^{t_0} \int (u_x^2+c_0 u^2)(t) \exp\left(\frac {\sqrt{c_0}} 4 (x-\rho(t_0)+\tfrac {c_0} 2 (t_0-t))\right)  dt	
\leq K.
\end{equation}
From \eqref{agarder}, there exists a sequence $t_n\to -\infty$ such that
$\int (u_x^2(t_n)+u^2(t_n)) \psi_0(t_n)\to 0$ as $n\to +\infty$.
Thus, from \eqref{futur2}--\eqref{futur3} applied to $t=t_n$, and passing to the limit as
$n\to +\infty$, we obtain
$$ \int_{-\infty}^{t_0} \int (u_{xx}^2 + c_0 u_x^2 + c_0^2 u^2)  \psi_{0x}    dt 
\leq   K \exp\left(-\frac {\sqrt{c_0}}4 {x_0}\right).
$$
Arguing as before with \eqref{truc}, we get, for all $t_0$,
\begin{equation}\label{agarder2}
\int_{-\infty}^{t_0} \int (u_{xx}^2+c_0u_x^2+c_0^2 u^2)(t) \exp\left(\frac {\sqrt{c_0}} 4 (x-\rho(t_0)+\tfrac {c_0} 2 (t_0-t))\right)  dt	
\leq K.
\end{equation}

Now, we use a monotonicity argument on $u_{xx}(t)$ to be able to give information on $v(t)$.
By similar calculations as in the proof of Lemma \ref{LE3as1}, we have
\begin{align*}
\label{}
   \frac d{dt} \int u_{xx}^2 \psi_0  &
   \leq - \int(3 u_{xxx}^2 + \tfrac {c_0} 2 u_{xx}^2 ) \psi_{0x} +
   \int u_{xx}^2 (f'(u)\psi_0 - f''(u) u_x \psi_0) \\&
   +2 \int u_{xx} u_x^2 f''(u) \psi_{0x} + 2 \int u_{xxx} u_x^2 f''(u) \psi_0.  
\end{align*}
We control the nonlinear terms as before, and then using \eqref{decayq},
$$
\int u_{xx}^2 (f'(u)\psi_0 - f''(u) u_x \psi_0)+2\int u_{xx} u_x^2 f''(u) \psi_{0x} \leq K \int (u_{xx}^2 + u_x^2) \psi_0,
$$
$$
\int u_{xxx} u_x^2 f''(u) \psi_0 \leq \int u_{xxx}^2 \psi_{0x} + K \int u_x^2 \psi_0.
$$
Arguing as before, we obtain the following conclusion, for all $t_0\in \mathbb{R}$,
\begin{equation}\label{prin20} \begin{split}
&\int (u_{xx}^2+ c_0 u_x^2+c_0^2 u^2)(t_0)\exp\left(\frac {\sqrt{c_0}} 4 (x-\rho(t_0))\right) 
\\& +
\int_{-\infty}^{t_0} \int (u_{xxx}^2+c_0 u_{xx}^2+ c_0^2u_x^2+c_0^3 u^2)(t) 
\exp\left(\frac {\sqrt{c_0}} 4 (x-\rho(t_0)+\tfrac {c_0} 2 (t_0-t))\right)  dt	
\leq K.
\end{split}\end{equation}
Since equation \eqref{fkdv} is invariant by the transformation $x\to -x$, $t\to -t$,  the
claim is proved.

\medskip

\noindent\emph{Proof of (\ref{monovNL}).}
Estimate \eqref{agarder4} and the decay on $Q_{c(t)}$ (see \eqref{decayQas1}) imply:
$$
\int v^2(t,x) \exp\left(\frac {\sqrt{c_0}} 4 |x|\right) dx \leq K.
$$
From this estimate, using the equation of $v$ and \eqref{prin20}, and arguing exactly as for the linear case
(proof of \eqref{monov}), we obtain  (\ref{monovNL}).


\begin{thebibliography}{10}

\bibitem{BL} H. Berestycki and P.-L. Lions,
Nonlinear scalar field equations. I. Existence of a ground state, Arch. Rational Mech. Anal. \textbf{ 82}, (1983) 313--345.
\bibitem{BSS} J.L. Bona, P.E. Souganidis and W.A. Strauss, Stability and instability of solitary waves of Korteweg--de Vries type, Proc. R. Soc. Lond. {\bf 411} (1987), 395--412.
\bibitem{BP}
V.S. Buslaev and  G. S. Perelman, On the stability of solitary waves for nonlinear Schršdinger equations.  Nonlinear evolution equations,  75--98, Amer. Math. Soc. Transl. Ser. 2, \textbf{164}, Amer. Math. Soc., Providence, RI, 1995.
\bibitem{ES} W. Eckhaus and P. Schuur, The emergence of solutions of the Korteweg--de Vries equation from arbitrary initial conditions,
Math. Meth. Appl. Sci., \textbf{5}, (1983) 97--116.
\bibitem{GRILLAKIS} M. Grillakis, Existence of nodal solutions of semilinear equations in $\mathbb{R}$, J. Diff. Eq. \textbf{85} (1990), 367--400.
\bibitem{KATO}  T. Kato, On the Cauchy problem for the (generalized) Korteweg-de Vries equation. Studies in applied mathematics, 93--128, Adv. Math. Suppl. Stud., \textbf{8}, Academic Press, New York, 1983.
\bibitem{KPV89} C.E. Kenig, G. Ponce and L. Vega,
On the (generalized) Korteweg-de Vries equation,
Duke Math. J. \textbf{59} (1989), 585--610.
\bibitem{KPV} C.E. Kenig, G. Ponce and L. Vega, Well-posedness and scattering results for the generalized Korteweg--de Vries equation via the contraction principle. Comm. Pure Appl. Math. \textbf{46} (1993), 527--620.
\bibitem{KRIEGERSCHLAG}
J. Krieger, W. Schlag, Stable manifolds for all monic supercritical focusing nonlinear Schršdinger equations in one dimension.  J. Amer. Math. Soc.  \textbf{19}  (2006),   815--920.
\bibitem{Martel} Y. Martel, Asymptotic $N$--soliton--like solutions of the subcritical and critical generalized Korteweg--de Vries equations, Amer. J. Math.  \textbf{127} (2005), 1103-1140.
\bibitem{yvanSIAM} Y. Martel, Linear problems related to asymptotic stability of solitons
of the generalized KdV equations, SIAM J. Math. Anal. \textbf{38} (2006), 759--781.
\bibitem{MMjmpa} Y. Martel and F. Merle, A Liouville theorem for the critical generalized Korteweg--de Vries equation. J. Math. Pures Appl. \textbf{79} (2000), 339--425.
\bibitem{MMarchives} Y. Martel and F. Merle, Asymptotic stability of solitons for subcritical generalized KdV equations. Arch. Ration. Mech. Anal. \textbf{157} (2001), 219--254.
\bibitem{MMgafa}  Y. Martel and F. Merle, Instability of solitons for the critical generalized Korteweg-de Vries equation. Geom. Funct. Anal. 11 (2001),  74--123.
\bibitem{MMnonlinearity} Y. Martel and F. Merle, Asymptotic stability of solitons of the subcritical gKdV equations revisited. Nonlinearity \textbf{18} (2005), no. 1, 55--80.
\bibitem{MMcol1} Y. Martel and F. Merle, 
Description of two soliton collision for the quartic KdV equation,
	arXiv:0709.2672v1.
\bibitem{MMcol2} Y. Martel and F. Merle, Stability of two soliton collision for the   gKdV equations, arXiv:0709.2677v1.
\bibitem{MMas2} Y. Martel and F. Merle, Refined asymptotics around solitons for the gKdV equations. 	To appear in Discrete and Continuous Dynamical Systems. Series A.
\bibitem{MMT} Y. Martel, F. Merle and Tai-Peng Tsai, Stability and asymptotic stability in the energy space of the sum of  $N$ solitons  for subcritical gKdV equations. Commun. Math. Phys. \textbf{231} (2002), 347--373.
\bibitem{MMTnls}
Y. Martel, F. Merle and  Tai-Peng Tsai, Stability in $H\sp 1$ of the sum of $K$ solitary waves for some nonlinear Schr\"odinger equations.  Duke Math. J.  \textbf{133}  (2006),   405--466.
\bibitem{Mizu} T. Mizumachi,
 Large time asymptotics of solutions around solitary waves to the generalized Korteweg-de Vries equations, SIAM J. Math. Anal. \textbf{32} (2001),  1050--1080.
\bibitem{PW} R.L. Pego and M.I. Weinstein, Asymptotic stability of solitary waves. Commun. Math. Phys. \textbf{164} (1994), 305--349.
\bibitem{P} G.S. Perelman, Asymptotic stability of multi-soliton solutions for nonlinear Schršdinger equations.  Comm. Partial Differential Equations \textbf{29}, (2004)  1051--1095.
\bibitem{RSS} I. Rodnianski, W. Schlag, A.D. Soffer, Asymptotic stability of $N$-soliton states of NLS, to appear in Comm. Pure. Appl. Math. 
\bibitem{SULEMSULEM} C. Sulem and P.-L. Sulem, The nonlinear Schr\"odinger equation.
Self-focusing and wave collapse. Applied Mathematical Sciences, \textbf{139}.
Springer-Verlag, New York, 1999.
\bibitem{We2} M.I. Weinstein, Modulational stability of ground states of nonlinear Schr\"odinger equations, SIAM J. Math. Anal. \textbf{16}, (1985) 472--491. 
\bibitem{We1} M.I. Weinstein, Lyapunov stability of ground states of nonlinear dispersive evolution equations. Comm. Pure Appl. Math. \textbf{39} (1986), 51--68.
\end{thebibliography}
\end{document}